# Estimating the joint distribution of independent categorical variables via model selection

C. DUROT[1], E. LEBARBIER[2] and A.-S. TOCQUET[3]

[1]*Laboratoire de mathématiques, Bât 425, Université Paris Sud, 91405 Orsay Cedex, France.
E-mail: Cecile.Durot@math.u-psud.fr*

[2]*Département MMIP, 16 rue Claude Bernard, 75231 Paris Cedex 05, France.
E-mail: emilie.lebarbier@agroparistech.fr*

[3]*Laboratoire Statistique et Génome, 523 place des Terrasses de l'Agora, F-91000 Evry, France.
E-mail: anne-sophie.tocquet@genopole.cnrs.fr*

Assume one observes independent categorical variables or, equivalently, one observes the corresponding multinomial variables. Estimating the distribution of the observed sequence amounts to estimating the expectation of the multinomial sequence. A new estimator for this mean is proposed that is nonparametric, non-asymptotic and implementable even for large sequences. It is a penalized least-squares estimator based on wavelets, with a penalization term inspired by papers of Birgé and Massart. The estimator is proved to satisfy an oracle inequality and to be adaptive in the minimax sense over a class of Besov bodies. The method is embedded in a general framework which allows us to recover also an existing method for segmentation. Beyond theoretical results, a simulation study is reported and an application on real data is provided.

## 1. Introduction

Assume that we observe independent categorical variables $Y_1, \ldots, Y_n$ which take values in a finite alphabet coded $\{1, \ldots, r\}$ with some integer $r \geq 2$. Such observations emerge, for instance, in DNA studies when the independence assumption is postulated. Equivalently, we observe the corresponding independent multinomial variables $X_1, \ldots, X_n$ where $X_i = (1, 0, \ldots, 0)^{\mathrm{T}}$ if $Y_i = 1$, and so forth. In this paper, we aim at estimating the distribution of the observed sequence. This amounts to estimating the expectation of the sequence $(X_1, \ldots, X_n)$ since the probability vector of $X_i$ is also its expectation. To make the connection with known methods from the literature, let us denote by $\mathcal{P}_r$ the set of all probabilities on $\{1,\ldots,r\}$, by $\mathcal{I}_n$ the set $\{i/n, i=1,\ldots,n\}$ and consider the function $\pi_n:\mathcal{I}_n \to \mathcal{P}_r$ such that for every $i$,

$$\pi_n(i/n) = E(X_i). \tag{1}$$

Our aim is to estimate $\pi_n$ under as little prior information as possible on this function. Note that $\pi_n$ depends on $n$ if $E(X_i)$ does not (which is typically the case for DNA sequences).

The paper in the literature whose issue is the closest to ours is the one by Aerts and Veraverbeke [1]. However, the approach in this paper is asymptotic; there, $Y_x$ denotes a categorical response variable at a fixed design point $x \in [0,1]^d$, and the aim is to estimate $\pi(x) = E(X_x)$ with $X_x$ the multinomial variable corresponding to $Y_x$ on the basis of independent observations $Y_{x_1},\ldots,Y_{x_n}$. The main difference between this model and ours is that $\pi$ does not depend on $n$, whence the approach is asymptotic. The similarity is that no specific model is postulated for this function. Assuming that $\pi$ is Lipschitz of order one, Aerts and Veraverbeke show that, at a fixed point $x \in (0,1)^d$, a kernel estimator of $\pi(x)$ with a proper bandwidth is asymptotically Gaussian as $n \to \infty$ and they propose a consistent Bootstrap procedure. As usual, the optimal choice for their bandwidth depends on the smoothness of $\pi$, the order $n^{-1/(4+d)}$ they propose for $d \leq 4$ being suitable only if $\pi$ is Lipschitz of order one.

From our knowledge, [1] is the only paper where no specific model is postulated. Most of the known nonparametric estimators in the literature for the joint distribution of multinomial variables are concerned with the segmentation model, which assumes that $\pi_n$ is a step function. Thus, there is a partition of $\{1,\ldots,n\}$ into segments within which the observations follow the same distribution and between which observations have different distributions, and an important task is to detect the change points. Within the segmentation model, several methods consist in minimizing, for every fixed $k$, a given contrast function over all possible partitions into $k$ segments, then selecting a convenient $k$ by penalizing the contrast function. Braun, Braun and Müller [8] consider a regression model that covers the case of multinomial variables. They assume $\pi_n$ independent of $n$ and penalize the quasi-deviance with a modified Schwarz criterion. Lebarbier and Nédélec [15] consider both least-squares and likelihood contrasts and propose a penalization inspired by Birgé and Massart [5], allowing $\pi_n$ to depend on $n$. In a slightly different spirit, Fu and Curnow [10] maximize the likelihood under restrictions on the lengths of the segments. The problem with these methods, in addition to being based on a quite restrictive model, is that they require a long computational time which does not scale to large sequences: using dynamic programming algorithms, the computational effort is $\mathcal{O}(n^2)$ when the number of change points is known, thus $\mathcal{O}(n^3)$ when this number is allowed to grow with $n$.

In this paper, we do not postulate a specific model for $\pi_n$. We propose an estimator which is adaptive, non-asymptotic, nonparametric and implementable even for large $n$. For every $i$ and $j$, we denote by $s_i^{(j)}$ the $j$th component of the vector $E(X_i)$ and we build a penalized contrast estimator as follows: Assume $n = 2^N$ for some integer $N \geq 1$ and let $\{\phi_\lambda, \lambda = 1,\ldots,n\}$ be the orthonormal Haar basis of $\mathbb{R}^n$. We consider a collection $\mathcal{M}$ of subsets of $\{1,\ldots,n\}$ and for every $m \in \mathcal{M}$, we compute the least-squares estimators of the vectors

$$s^{(j)} := (s_1^{(j)},\ldots,s_n^{(j)}), \qquad j = 1,\ldots,r$$

under the assumption that they all belong to the linear span of $\{\phi_\lambda, \lambda \in m\}$. Then we select a convenient subset $m$ by penalizing the least-squares contrast after the fashion of Birgé and Massart [5]. At first, we consider an exhaustive strategy which roughly consists in defining $\mathcal{M}$ as the collection of all possible subsets of $\{1,\ldots,n\}$. We prove that, up to a $\log n$ factor, the resulting estimator satisfies an oracle inequality and is adaptive in the minimax sense over a class of Besov bodies. Subsequently, we consider a non-exhaustive strategy: Based on the compression algorithm developed in [4], we build a collection $\mathcal{M}$ in such a way that the resulting estimator indeed satisfies an oracle inequality and is adaptive over this class of Besov bodies (without $\log n$ factor). Both strategies are implementable even for large $n$ since all the least-squares estimators can be computed at the same time: The computational effort is $\mathcal{O}(n \log n)$ with both strategies.

In order to gain generality, we embed our method into a general framework: We consider a collection $\{S_m, m \in \mathcal{M}\}$ of linear subspaces of $\mathbb{R}^n$, we estimate the vectors $s^{(j)}$, $j = 1,\ldots,r$ under the assumption that they all belong to $S_m$ and then select a convenient model $S_m$. Thus for the above strategies, $S_m$ is merely the span of $\{\phi_\lambda, \lambda \in m\}$ for a given $m \subset \{1,\ldots,n\}$. The general framework also recovers the method of Lebarbier and Nédélec [15] and allows us to slightly improve some of their results. Moreover, we prove that their method is adaptive in the segmentation model.

The paper is organized as follows: The general framework is described in Section 2. Subsequently, special cases are more precisely described: the method of Lebarbier and Nédélec is studied in Section 3 and the exhaustive and non-exhaustive strategies based on the Haar basis are studied in Sections 4 and 5, respectively. A simulation study is reported in Section 6 and an application on a DNA sequence is given in Section 7. The proof of the main results is postponed to Section 8.

## 2. The general framework

Assume one observes independent variables $X_1,\ldots,X_n$ where $X_i \in \{0,1\}^r$ with some integer $r \geq 2$ and the components of $X_i$ sum up to one. Hence if one observes independent categorical variables $Y_1,\ldots,Y_n$ which take values in $\{1,\ldots,r\}$, the $j$th component of $X_i$ equals one if $Y_i = j$ and zero otherwise. We aim at estimating the expectation of the sequence $(X_1,\ldots,X_n)$. Notation given in Section 2.1 will be used throughout the paper. Our penalized least-squares estimator is defined in Section 2.2 with an arbitrary penalty function. The choice of the form of the penalty function is discussed in Section 2.3.

### 2.1. Notation

For our statistical model, we use the following notation: Let $X$ be the $r \times n$ random matrix with $i$th column $X_i$, let $s$ be the expectation of $X$ and set $\varepsilon = X - s$. The underlying probability and expectation in our model are denoted by $\mathbb{P}_s$ and $\mathbb{E}_s$, respectively.

For vectors and matrices, we use the following notation: By convention, vectors in $\mathbb{R}^r$ are column vectors while vectors in $\mathbb{R}^n$ are line vectors. The $j$th line of a matrix $t$ is denoted by $t^{(j)}$ while its $i$th column is denoted by $t_i$. The Euclidean inner product and Euclidean norm in $\mathbb{R}^p$, $p \in \mathbb{N}$, are denoted by $\langle \cdot, \cdot \rangle_p$ and $\| \cdot \|_p$, respectively. The space of $r \times n$ matrices is equipped with the inner product

$$\langle t, t' \rangle = \sum_{i,j} t_i^{(j)} t'^{(j)}_i$$

and we denote by $\| \cdot \|$ the associated norm. For every linear subspace $S$ of $\mathbb{R}^n$, $\mathbb{R}^r \otimes S$ denotes the set of those $r \times n$ matrices the lines of which belong to $S$. Finally, we use the same notation for a vector in $\mathbb{R}^n$ and for the corresponding discrete function. Thus, $v$ may denote the vector $(v_1, \ldots, v_n)$ as well as the discrete function $v : i \mapsto v_i$, $i = 1, \ldots, n$.

## 2.2. The estimator

Recall that we aim at estimating $s$ on the basis of $X_1, \ldots, X_n$. For this task, consider a finite collection $\{S_m, m \in \mathcal{M}\}$ of linear subspaces of $\mathbb{R}^n$ where $\mathcal{M}$ is allowed to depend on $n$. For every $m \in \mathcal{M}$, let $\hat{s}_m$ be the least-squares estimator of $s$ under the assumption that the lines of $s$ all belong to $S_m$: $\hat{s}_m$ is the unique minimizer of $\|t - X\|^2$ over $t \in \mathbb{R}^r \otimes S_m$, that is, the orthogonal projection of $X$ onto $\mathbb{R}^r \otimes S_m$. For ease of computation, it is worth noticing that

$$\|t - X\|^2 = \sum_{j=1}^{r} \|t^{(j)} - X^{(j)}\|_n^2,$$

hence for every $m$ and $j$, $\hat{s}_m^{(j)}$ is the orthogonal projection of $X^{(j)}$ onto $S_m$. It is also worth noticing that $\hat{s}_m$ minimizes the contrast

$$\gamma_n(t) = \|t\|^2 - 2\langle X, t \rangle$$

over $t \in \mathbb{R}^r \otimes S_m$. The estimator we are interested in is the penalized contrast estimator defined by

$$\tilde{s} = \hat{s}_{\hat{m}},$$

where $\hat{m}$ minimizes $\gamma_n(\hat{s}_m) + \text{pen}(m)$ over $m \in \mathcal{M}$. The function $\text{pen} : \mathcal{M} \to \mathbb{R}_+$ is called the penalty function and remains to be chosen, see Section 2.3. Typically, $\text{pen}(m)$ increases with the dimension of $S_m$. The performances of our estimator then depend on the strategy, that is, the way we choose the collection of models and the penalty function. Strategies where $S_m$ is generated by either step functions or Haar functions are studied in detail in subsequent sections.

Let us notice that since we do not take into account the constraint that every column of $s$ is a probability on $\{1, \ldots, r\}$, our estimator may not satisfy this constraint: The components of a given column of the estimator sum up to one provided the vector $(1, \ldots, 1)$ belongs to each $S_m$, $m \in \mathcal{M}$, but it may happen that components of the estimator do not belong to $[0,1]$. However one may always project $\tilde{s}$ on the closed convex subset of those matrices that satisfy the constraint, getting an estimator that is closer to $s$ than $\tilde{s}$. Note also that the strategies studied in Sections 3–5 are such that $(1,\ldots,1)$ indeed belongs to each space $S_m$, $m \in \mathcal{M}$. Thus, we assume in Section 2.3 that the dimension of $S_m$ is $D_m \geq 1$ for each $m$.

## 2.3. The form of the penalty function

The following theorem provides an upper bound for the $L_2$-risk of our estimator, which gives indications on how to choose the penalty function.

**Theorem 1.** *Let $\{S_m, m \in \mathcal{M}\}$ be a finite collection of linear subspaces of $\mathbb{R}^n$, for every $m \in \mathcal{M}$, let $D_m$ be the dimension of $S_m$ and assume $D_m \geq 1$. Let* $\mathrm{pen}\colon \mathcal{M} \to \mathbb{R}_+$ *satisfy*

$$\mathrm{pen}(m) \geq K D_m (1 + 2\sqrt{2L_m})^2 \qquad \forall m \in \mathcal{M} \tag{2}$$

*for some $K > 1$ and some positive numbers $L_m$. Then*

$$\mathbb{E}_s \|s - \tilde{s}\|^2 \leq \frac{2K(K+1)^2}{(K-1)^3} \inf_{m \in \mathcal{M}} \{\|s - \bar{s}_m\|^2 + \mathrm{pen}(m)\} + 16K\left(\frac{K+1}{K-1}\right)^3 \Sigma,$$

*where $\bar{s}_m$ is the orthogonal projection of $s$ onto $\mathbb{R}^r \otimes S_m$ and*

$$\Sigma = \sum_{m \in \mathcal{M}} \exp(-L_m D_m). \tag{3}$$

This suggests a choice for the form of the penalty function: It can be chosen in such a way that the upper bound of the risk is as small as possible, subject to the constraint (2). A possible choice is

$$\mathrm{pen}(m) = D_m(k_1 L_m + k_2) \tag{4}$$

for some positive $k_1$, $k_2$ and $L_m$, where the weights $L_m$ have to be chosen as small as possible subject to the constraint that $\Sigma$ be small, $\Sigma \leq 1$, say. In order to check the relevance of this choice, we aim at proving an oracle-type inequality: It is expected that the risk of $\tilde{s}$ is close to the minimal risk of the least-squares estimators $\hat{s}_m$, in the sense that

$$\mathbb{E}_s \|s - \tilde{s}\|^2 \leq C \inf_{m \in \mathcal{M}} \mathbb{E}_s \|s - \hat{s}_m\|^2, \tag{5}$$

where, ideally, $C > 0$ does not depend on $r$ and $n$.

**Corollary 1.** *Let $\{S_m, m \in \mathcal{M}\}$ be a finite collection of linear subspaces of $\mathbb{R}^n$, for every $m \in \mathcal{M}$, let $D_m$ be the dimension of $S_m$ and assume $D_m \geq 1$. Let* $\mathrm{pen}\colon \mathcal{M} \to \mathbb{R}_+$ *be defined by (4) with some positive $k_1$, $k_2$ and $L_m$'s that satisfy $\Sigma \leq 1$, see (3). Assume $s_i^{(j)} \leq 1 - \rho$* for some $\rho > 0$ and all $i$, $j$. If $k_1$ and $k_2$ are large enough, then there exists $C > 0$ that only depends on $k_1$, $k_2$ and $\rho$ such that

$$\mathbb{E}_s \|s - \tilde{s}\|^2 \leq C \left(1 + \sup_{m \in \mathcal{M}} L_m \right) \inf_{m \in \mathcal{M}} \mathbb{E}_s \|s - \hat{s}_m\|^2.$$

Thus (5) holds with $C$ depending on neither $n$ nor $r$, provided $k_1$, $k_2$ are large enough and $L_m \leq L$ for every $m$ and an absolute constant $L$. However, $\mathcal{M}$ is allowed to depend on $n$ and it is assumed that $\Sigma \leq 1$ so this could be possible only if the collection of models is not too rich: If $\mathcal{M}$ contains many models with dimension $D_m$, then $L_m$ has to increase with $n$ to ensure $\Sigma \leq 1$ and (5) only holds with $C > 0$ depending on $n$ in that case.

**Proof of Corollary 1.** For every $m$, let $\{\psi_{m1}, \ldots, \psi_{mD_m}\}$ be an orthonormal basis of $S_m$. For every $j$, the $j$th line of $\hat{s}_m$ is the orthogonal projection of $X^{(j)}$ onto $S_m$ hence

$$\hat{s}_m^{(j)} = \sum_{k=1}^{D_m} \langle X^{(j)}, \psi_{mk} \rangle_n \psi_{mk}. \tag{6}$$

Since $\bar{s}_m$ is defined in the same way as $\hat{s}_m$ with $X$ replaced by $s$, we get

$$\|\bar{s}_m - \hat{s}_m\|^2 = \sum_{j=1}^{r} \sum_{k=1}^{D_m} (\langle \varepsilon^{(j)}, \psi_{mk} \rangle_n)^2.$$

We have $\text{var}_s(\varepsilon_i^{(j)}) = s_i^{(j)}(1 - s_i^{(j)})$, where $1 - s_i^{(j)} \geq \rho$ and the components of $s_i$ sum up to one so

$$\sum_{j=1}^{r} \text{var}_s(\varepsilon_i^{(j)}) \geq \rho.$$

Moreover, the components of $\varepsilon^{(j)}$ are independent and centered, hence

$$\mathbb{E}_s \|\bar{s}_m - \hat{s}_m\|^2 \geq \rho D_m.$$

It thus follows from Pythagoras' equality that

$$\mathbb{E}_s \|s - \hat{s}_m\|^2 \geq \|s - \bar{s}_m\|^2 + \rho D_m.$$

By Theorem 1, there exists $C' > 0$ that only depends on $k_1$ and $k_2$ such that

$$\mathbb{E}_s \|s - \tilde{s}\|^2 \leq C' \left(1 + \sup_{m \in \mathcal{M}} L_m \right) \inf_{m \in \mathcal{M}} \{\|s - \bar{s}_m\|^2 + D_m\}.$$

The result easily follows from the last two displays. □

## 3. Exhaustive strategy based on indicator functions

In the segmentation model, it is assumed that there exists a partition of $\{1,\ldots,n\}$ into intervals within which the observations follow the same distribution and between which observations have different distributions. This amounts to assuming that the functions $s^{(j)}$, $j=1,\ldots,r$ are piecewise-constant with respect to a common partition of $\{1,\ldots,n\}$, so one can use in this model the general method of Section 2, as follows.

For every subset $I \subset \{1,\ldots,n\}$, let $\psi_I$ denote the indicator function of $I$, which means that $\psi_I$ is a vector in $\mathbb{R}^n$ with an $i$th component equal to one if $i \in I$ and zero otherwise. Let $\mathcal{M}$ be the collection of all partitions of $\{1,\ldots,n\}$ into intervals and for every $m \in \mathcal{M}$ let $S_m$ be the linear span of $\{\psi_I, I \in m\}$. The dimension of $S_m$ is the cardinality of $m$, hence for every $D=1,\ldots,n$ the number $N_D$ of models with dimension $D$ satisfies

$$N_D = \binom{n-1}{D-1} \leq \binom{n}{D} \leq \left(\frac{en}{D}\right)^D. \tag{7}$$

Setting $L_m = 2 + \log(n/D_m)$ thus ensures $\Sigma \leq 1$ (see (3)):

$$\sum_{m\in\mathcal{M}} \exp(-L_m D_m) = \sum_{D=1}^{n} N_D \exp\left(-\left(2+\log\left(\frac{n}{D}\right)\right)D\right) \leq 1. \tag{8}$$

By Corollary 1, the estimator based on a penalty function of the form (4) thus satisfies an oracle inequality up to a $\log n$ factor, provided $k_1$ and $k_2$ are large enough. With the above definition of $L_m$, the penalty function takes the form

$$\text{pen}(m) = D_m\left(c_1 \log\left(\frac{n}{D_m}\right) + c_2\right). \tag{9}$$

In the sequel, we denote by EI (an abbreviation of Exhaustive/Indicator) the strategy based on the above collection $\mathcal{M}$ and a penalty function of the form (9). The following proposition is then easily derived from Corollary 1.

**Proposition 1.** *Let $\tilde{s}$ be computed with the EI strategy. Assume $s_i^{(j)} \leq 1-\rho$, $\rho > 0$, for all $i,j$. If $c_1$, $c_2$ are large enough, then there exists $C>0$ that only depends on $\rho$, $c_1$, $c_2$ such that*

$$\mathbb{E}_s \|s - \tilde{s}\|^2 \leq C \log n \inf_{m\in\mathcal{M}} \mathbb{E}_s \|s - \hat{s}_m\|^2.$$

The EI strategy was first studied by Lebarbier and Nédélec [15]. Proposition 1 is a slight improvement of one of their results since they obtain this inequality with $C$ replaced by $Cr$ and under a more restrictive assumption. With this strategy, for every $I \in \hat{m}$ and every $i \in I$, the $i$th column of $\tilde{s}$ is the mean of the vectors $(X_l)_{l\in I}$ so, in particular, every column of $\tilde{s}$ is a probability on $\{1,\ldots,r\}$.

We pursue the study of the EI strategy by proving that it is adaptive in the segmentation model. For a given positive integer $k$, let $\mathcal{P}_k$ be the set of those $r \times n$ matrices $t$ such that every column of $t$ is a probability on $\{1,\ldots,r\}$ and there exists a partition of $\{1,\ldots,n\}$ into $k$ intervals on which the functions $t^{(j)}$, $j=1,\ldots,r$, are constant. In the sequel, we prove that the EI strategy is adaptive with respect to $k$ in the sense that the corresponding estimator is minimax simultaneously over the sets $\mathcal{P}_k$, $k=3,\ldots,n$. A lower bound for the minimax risk over $\mathcal{P}_k$ is given in the following theorem and the minimax result is derived in Corollary 2.

**Theorem 2.** *Assume $n \geq 3$ and fix $k \in \{3,\ldots,n\}$. There exists an absolute constant $C > 0$ such that*

$$\inf_{\hat{s}} \sup_{s \in \mathcal{P}_k} \mathbb{E}_s \|s - \hat{s}\|^2 \geq Ck\left(\log\left(\frac{n}{k}\right) + 1\right),$$

*where the infimum is taken over all estimators $\hat{s}$.*

**Corollary 2.** *Let $\tilde{s}$ be computed with the EI strategy. If $c_1$, $c_2$ are large enough, $n \geq 3$ and $k \in \{3,\ldots,n\}$, then there exists $C > 0$ that only depends on $c_1$ and $c_2$ such that*

$$\sup_{s \in \mathcal{P}_k} \mathbb{E}_s \|s - \tilde{s}\|^2 \leq C \inf_{\hat{s}} \sup_{s \in \mathcal{P}_k} \mathbb{E}_s \|s - \hat{s}\|^2.$$

**Proof.** We may assume $\Sigma \leq 1$ so Theorem 1 shows that there exists $C > 0$ that only depends on $c_1$ and $c_2$ such that for every $s$

$$\mathbb{E}_s \|s - \tilde{s}\|^2 \leq C \min_D \left\{ \min_{m \in \mathcal{M}, D_m = D} \|s - \bar{s}_m\|^2 + D\left(\log\left(\frac{n}{D}\right) + 1\right) \right\}$$

$$\leq C \left\{ \min_{m \in \mathcal{M}, D_m = k} \|s - \bar{s}_m\|^2 + k\left(\log\left(\frac{n}{k}\right) + 1\right) \right\}.$$

The latter minimum vanishes if $s \in \mathcal{P}_k$, so the result follows from Theorem 2.　□

## 4. Exhaustive strategy based on Haar functions

In this section, we assume $n = 2^N$ for an integer $N \geq 1$ and we consider models generated by Haar functions, as defined below.

**Definition 1.** *Let $\Lambda = \bigcup_{j=-1}^{N-1} \Lambda(j)$, where $\Lambda(-1) = \{(-1,0)\}$ and*

$$\Lambda(j) = \{(j,k), k = 0,\ldots, 2^j - 1\}$$

*for $j \neq -1$. Let $\varphi: \mathbb{R} \to \{-1,1\}$ be the function with support $(0,1]$ that takes value 1 on $(0,1/2]$ and $-1$ on $(1/2,1]$. For every $\lambda \in \Lambda$, let $\phi_\lambda$ be the vector in $\mathbb{R}^n$ with $i$th component $\phi_{\lambda i} = 1/\sqrt{n}$ if $\lambda = (-1,0)$ and*

$$\phi_{\lambda i} = \frac{2^{j/2}}{\sqrt{n}} \varphi\left(2^j \frac{i}{n} - k\right)$$

*if $\lambda = (j,k)$ for some $j \neq -1$. The orthonormal basis $\{\phi_\lambda, \lambda \in \Lambda\}$ of $\mathbb{R}^n$ is called the Haar basis and its elements are called Haar functions.*

Let $\mathcal{M}$ be the collection of all subsets of $\Lambda$ that contain $(-1,0)$ and for every $m \in \mathcal{M}$, let $S_m$ be the linear span of $\{\phi_\lambda, \lambda \in m\}$. The number $N_D$ of models with dimension $D$ satisfies (7). In a manner similar to that of Section 3, we denote by EH (an abbreviation of Exhaustive/Haar) the strategy based on the above collection $\mathcal{M}$ and a penalty function of the form (9). The following proposition is then easily derived from Corollary 1.

**Proposition 2.** *Let $\tilde{s}$ be computed with the EH strategy. Assume $s_i^{(j)} \leq 1 - \rho$, $\rho > 0$, for all $i, j$. If $c_1$, $c_2$ are large enough, then there exists $C > 0$ that only depends on $\rho$, $c_1$, $c_2$ such that*

$$\mathbb{E}_s \|s - \tilde{s}\|^2 \leq C \log n \inf_{m \in \mathcal{M}} \mathbb{E}_s \|s - \hat{s}_m\|^2.$$

A great advantage of the EH strategy is that the obtained estimator is easily computable in practice, as we explain now. For every $m$, let $\hat{s}_m$ be the orthogonal projection of $X$ onto $\mathbb{R}^r \otimes S_m$. Then $\gamma_n(\hat{s}_m) = -\|\hat{s}_m\|^2$ and it follows from (6) that $\|\hat{s}_m\|^2 = \sum_{\lambda \in m} \|\hat{\beta}_\lambda\|_r^2$, where

$$\|\hat{\beta}_\lambda\|_r^2 = \sum_{j=1}^r (\hat{\beta}_\lambda^{(j)})^2 \quad \text{and} \quad \hat{\beta}_\lambda^{(j)} = \langle X^{(j)}, \phi_\lambda \rangle_n.$$

For every $m$, $\text{pen}(m)$ only depends on $D_m$ so we can write $\text{pen}(m) = \text{pen}'(D_m)$. Hence

$$\min_{m \in \mathcal{M}} \{\gamma_n(\hat{s}_m) + \text{pen}(m)\} = \min_{1 \leq D \leq n} \left\{ -\max_{m, D_m = D} \sum_{\lambda \in m} \|\hat{\beta}_\lambda\|_r^2 + \text{pen}'(D) \right\}$$

$$= \min_{1 \leq D \leq n} \left\{ -\sum_{k=1}^D \|\hat{\beta}_{\tau(k)}\|_r^2 + \text{pen}'(D) \right\},$$

where $\|\hat{\beta}_{\tau(1)}\|_r^2 \geq \cdots \geq \|\hat{\beta}_{\tau(n)}\|_r^2$ are the coefficients arranged in decreasing order. The computation of $\tilde{s}$ thus requires a computational effort of order $\mathcal{O}(n \log n)$ using the following algorithm.

Algorithm EH: – Compute $\hat{\beta}_\lambda^{(j)}$ for every $j$ and $\lambda \in \Lambda$.
- Compute the coefficients $\|\hat{\beta}_\lambda\|_r^2$ for every $\lambda \in \Lambda$.
- Arrange these coefficients in decreasing order.
- Compute $-\sum_{k=1}^{D} \|\hat{\beta}_{\tau(k)}\|_r^2 + \mathrm{pen}'(D)$ for every $D = 1, \ldots, n$.
- Find $\hat{D}$ that minimizes the quantity above.
- Compute $\tilde{s}^{(j)} = \sum_{k=1}^{\hat{D}} \hat{\beta}_{\tau(k)}^{(j)} \phi_{\tau(k)}$.

We establish below that the EH strategy is adaptive up to a $\log n$ factor with respect to classes $\mathcal{P}(R, \alpha, p)$ of matrices related to Besov bodies, see Definition 2. To motivate the choice of these classes, we make the connection with the more usual problem of estimating a function $\pi \in L_2[0,1]$ on the basis of the observations

$$y_i = \pi(i/n) + \delta_i, \qquad i = 1, \ldots, n, \tag{10}$$

where the $\delta_i$'s are i.i.d. centered random variables. With the notation of Definition 1, define $\Phi_\lambda(x) = 2^{j/2} \varphi(2^j x - k)$ for all $j \geq 0$ and $\lambda = (j, k) \in \Lambda(j)$, and let $\Phi_{(-1,0)}(x) = 1$, $x \in [0,1]$. Then, $\{\Phi_\lambda, \lambda \in \cup_{j \geq -1} \Lambda(j)\}$ is the Haar basis of $L_2[0,1]$ and the coefficients of $\pi$ in this basis are

$$\tau_\lambda = \langle \pi, \Phi_\lambda \rangle_{L_2}$$

with $\langle \cdot, \cdot \rangle_{L_2}$ the scalar product in $L_2[0,1]$. If $\pi$ belongs to a Besov ball in $B_p^\alpha(L_p)$ with $\alpha < 1$, $\alpha > 1/p - 1/2$ and $p \leq 2$, then typically,

$$\sum_{j \geq 0} 2^{jp(\alpha + 1/2 - 1/p)} \sum_{\lambda \in \Lambda(j)} |\tau_\lambda|^p \leq R^p \tag{11}$$

for some $R > 0$, see Donoho and Johnstone [9] or Section 4.5 of Massart [16]. Constraints like (11) are thus often considered to derive minimax rates of convergence over Besov balls in model (10). Now, suppose there is $\pi \in L_2[0,1]$ such that for every $i$ and a fixed $j$,

$$s_i^{(j)} = \pi(i/n).$$

Then with the notation of Definition 1,

$$\beta_\lambda^{(j)} := \langle s^{(j)}, \phi_\lambda \rangle_n = \frac{1}{\sqrt{n}} \sum_{i=1}^{n} \pi(i/n) \Phi_\lambda(i/n).$$

Approximating an integral over $[0,1]$ with a discrete sum thus yields

$$\beta_\lambda^{(j)} \approx \sqrt{n} \tau_\lambda. \tag{12}$$

Thus, by analogy with (11), we consider the following classes $\mathcal{P}(R, \alpha, p)$.

**Definition 2.** *Let $R > 0$, $p \in (0,2]$, $\alpha > 1/p - 1/2$. Let $\mathcal{B}(R, \alpha, p)$ be the set of matrices $t \in \mathbb{R}^r \otimes \mathbb{R}^n$ such that the coefficients $\beta_\lambda^{(j)} = \langle t^{(j)}, \phi_\lambda \rangle_n$, $j = 1, \ldots, r$ satisfy*

$$\sum_{j=0}^{N-1} 2^{jp(\alpha+1/2-1/p)} \sum_{\lambda \in \Lambda(j)} \|\beta_\lambda\|_r^p \leq n^{p/2} R^p.$$

*We define $\mathcal{P}(R, \alpha, p)$ to be the set of those $r \times n$ matrices $s \in \mathcal{B}(R, \alpha, p)$ such that every column of $s$ is a probability on $\{1, \ldots, r\}$.*

We aim to prove that the EH strategy is adaptive (up to a $\log n$ factor) in the sense that the corresponding estimator is minimax (up to a $\log n$ factor) simultaneously over sets $\mathcal{P}(R, \alpha, p)$ for various $R$, $\alpha$ and $p$. In this regard, we first establish a lower bound for the minimax risk over these sets.

**Theorem 3.** *Let $R > 0$, $p \in (0, 2]$ and $\alpha > 1/p - 1/2$. Assume $n^\alpha > R$ and $nR^2 \geq 1$. Then there exists an absolute constant $C > 0$ such that*

$$\inf_{\hat{s}} \sup_{s \in \mathcal{P}(R,\alpha,p)} \mathbb{E}_s \|s - \hat{s}\|^2 \geq C(nR^2)^{1/(2\alpha+1)},$$

*where the infimum is taken over all estimators $\hat{s}$.*

In the problem of estimating $\pi \in L_2[0,1]$ within model (10), the minimax rate of convergence in the squared $L_2$-norm over sets like (11) is typically

$$R^{2/(2\alpha+1)} n^{-2\alpha/(2\alpha+1)}, \tag{13}$$

see Donoho and Johnstone [9] or Section 4.3 of Massart [16]. Similar to (12), $\|s^{(j)} - \hat{s}^{(j)}\|_n^2$ is approximately $n\|\pi - \hat{\pi}\|_{L_2}^2$ with $\hat{\pi}$ as an estimator of $\pi$ and $\|\cdot\|_{L_2}$ as the $L_2$-norm, so the lower bound in Theorem 3 compares to (13) and appears natural.

An upper bound for the maximal risk of our estimator is given in Theorem 4 below and the minimax result is derived in Corollary 3.

**Theorem 4.** *Let $\tilde{s}$ be computed with the EH strategy. Let $R > 0$, $p \in (0, 2]$, $\alpha > 1/p - 1/2$ and let $\rho$ be the unique solution in $(1, \infty)$ of $\rho \log \rho = 2\alpha + 1$. Assume $c_1$, $c_2$ large enough, $nR^2 \geq \log n > 1$ and $n^\alpha > R\sqrt{\rho}$. Then there exists $C > 0$ that only depends on $\alpha$, $p$, $c_1$ and $c_2$ such that*

$$\sup_{s \in \mathcal{P}(R,\alpha,p)} \mathbb{E}_s \|s - \tilde{s}\|^2 \leq C(nR^2)^{1/(2\alpha+1)} \left[1 + \log\left(\frac{n^\alpha}{R}\right)\right]^{2\alpha/(2\alpha+1)}.$$

**Corollary 3.** *Under the assumptions of Theorem 4, there exists $C > 0$ that only depends on $\alpha$, $p$, $c_1$ and $c_2$ such that*

$$\sup_{s \in \mathcal{P}(R,\alpha,p)} \mathbb{E}_s \|s - \tilde{s}\|^2 \leq C \left[1 + \log\left(\frac{n^\alpha}{R}\right)\right]^{2\alpha/(2\alpha+1)} \inf_{\hat{s}} \sup_{s \in \mathcal{P}(R,\alpha,p)} \mathbb{E}_s \|s - \hat{s}\|^2.$$

## 5. A non-exhaustive strategy based on Haar functions

The strategy described in the preceding section is attractive because it is easily implementable. It satisfies an oracle inequality and is minimax over Besov bodies up to a log factor. In this section, we propose another strategy based on Haar functions that is also easily implementable and indeed satisfies an oracle inequality and a minimax result over Besov bodies (without log factor). With this strategy, $\mathcal{M}$ is composed of some well-chosen subsets of $\Lambda$ that contain $(-1,0)$, whereas the exhaustive strategy of the preceding section considers all possible such subsets. As a consequence, the penalty function can be chosen smaller than in the preceding section which allows us to select models with higher dimensions. The obtained estimator thus should have better risk than in the preceding section, provided the proposed collection of models properly approximates the exhaustive collection. For the choice of a proper collection of models, we were inspired by Birgé and Massart ([4] and Section 6.4 of [5]).

In this section, we assume $n = 2^N$ for an integer $N \geq 1$ and we use the same notation as in Definition 1. For every $J \in \{0, \ldots, N-1\}$ let

$$\mathcal{M}_J = \left\{ m \subset \Lambda, m = \left[ \bigcup_{j=-1}^{J-1} \Lambda(j) \right] \cup \left[ \bigcup_{k=0}^{N-J-1} \Lambda'(J+k) \right] \right\},$$

where $\Lambda'(J+k)$ is any subset of $\Lambda(J+k)$ with cardinality $\lfloor 2^J(k+1)^{-3} \rfloor$ (here, $\lfloor x \rfloor$ is the integer part of $x$). Let

$$\mathcal{M} = \bigcup_{J=0}^{N-1} \mathcal{M}_J$$

and, for every $m \in \mathcal{M}$, let $S_m$ be the linear span of $\{\phi_\lambda, \lambda \in m\}$. We have $D_m \geq 2^J$ for every $m \in \mathcal{M}_J$ and, from Proposition 4 of [4], there exists an absolute constant $A$ such that the cardinality of $\mathcal{M}_J$ is less than $\exp(A2^J)$. Setting $L_m = L$ for every $m \in \mathcal{M}$ and some large enough absolute constant $L$ thus ensures $\Sigma \leq 1$, see (3). With this definition of $L_m$, a penalty function of the form (4) takes the form

$$\text{pen}(m) = cD_m \tag{14}$$

for some positive number $c$. In the sequel, we denote by NEH (an abbreviation of non-exhaustive/Haar) the strategy based on the above collection $\mathcal{M}$ and a penalty function of the form (14). It follows from Corollary 1 that the corresponding estimator satisfies an oracle inequality.

**Proposition 3.** *Let $\tilde{s}$ be computed with the NEH strategy. Assume $s_i^{(j)} \leq 1 - \rho$, $\rho > 0$, for all $i, j$. If $c$ is large enough, then there exists $C > 0$ that only depends on $\rho$ and $c$ such that*

$$\mathbb{E}_s \|s - \tilde{s}\|^2 \leq C \inf_{m \in \mathcal{M}} \mathbb{E}_s \|s - \hat{s}_m\|^2.$$

The NEH strategy is easily implementable, as we explain now. All elements of a given $\mathcal{M}_J$ have the same cardinality so for every $m \in \mathcal{M}_J$, $\mathrm{pen}(m)$ only depends on $J$. Similar to Section 4, we can write $\mathrm{pen}(m) = \mathrm{pen}'(J)$, hence

$$\min_{m \in \mathcal{M}} \{\gamma_n(\hat{s}_m) + \mathrm{pen}(m)\} = \min_{0 \leq J \leq N-1} \left\{ -\max_{m \in \mathcal{M}_J} \sum_{\lambda \in m} \|\hat{\beta}_\lambda\|_r^2 + \mathrm{pen}'(J) \right\}$$

$$= \min_{0 \leq J \leq N-1} \left\{ -\sum_{\lambda \in \hat{m}_J} \|\hat{\beta}_\lambda\|_r^2 + \mathrm{pen}'(J) \right\}.$$

Here,

$$\hat{m}_J = \left[ \bigcup_{j=-1}^{J-1} \Lambda(j) \right] \cup \left[ \bigcup_{k=0}^{N-J-1} \tilde{\Lambda}(J+k) \right],$$

where $\tilde{\Lambda}(J+k)$ is the set that contains the indices of the $\lfloor 2^J(k+1)^{-3} \rfloor$ bigger coefficients $\|\hat{\beta}_\lambda\|_r$, $\lambda \in \Lambda(J+k)$. The computation of $\tilde{s}$ thus requires a computational effort of order $\mathcal{O}(n \log n)$ using the following algorithm.

Algorithm NEH: – Compute $\hat{\beta}_\lambda^{(j)}$ for every $j$ and $\lambda \in \Lambda$.
– Compute the coefficients $\|\hat{\beta}_\lambda\|_r^2$ for every $\lambda \in \Lambda$.
– For every $j$, arrange $\|\hat{\beta}_\lambda\|_r, \lambda \in \Lambda(j)$ in decreasing order.
– Compute $\hat{m}_J$ for every $J$.
– Compute $-\sum_{\lambda \in \hat{m}_J} \|\hat{\beta}_\lambda\|_r^2 + \mathrm{pen}'(J)$ for every $J$.
– Find $\hat{J}$ that minimizes the quantity above.
– Compute $\tilde{s}^{(j)} = \sum_{\lambda \in \hat{m}_{\hat{J}}} \hat{\beta}_\lambda^{(j)} \phi_\lambda$.

To conclude the study of the NEH strategy, we prove that it is adaptive in the sense that the corresponding estimator is minimax simultaneously over sets $\mathcal{P}(R, \alpha, p)$ for various $R$, $\alpha$ and $p$, see Definition 2. The following theorem provides an upper bound for the maximal risk of our estimator. Comparing this upper bound with the lower bound in Theorem 3 proves the minimax result stated in Corollary 4.

**Theorem 5.** *Let $\tilde{s}$ be computed with the NEH strategy. Let $R > 0$, $p \in (0, 2]$ and $\alpha > 1/p - 1/2$. Assume $c$ large enough, $n^\alpha > R$ and $nR^2 \geq 1$. Then there exists $C > 0$ that only depends on $\alpha$, $p$ and $c$ such that*

$$\sup_{s \in \mathcal{P}(R, \alpha, p)} \mathbb{E}_s \|s - \tilde{s}\|^2 \leq C(nR^2)^{1/(2\alpha+1)}.$$

**Corollary 4.** *Assume the assumptions of Theorem 5 are fulfilled. Then there exists $C > 0$ that only depends on $\alpha$, $p$ and $c$ such that*

$$\sup_{s \in \mathcal{P}(R, \alpha, p)} \mathbb{E}_s \|s - \tilde{s}\|^2 \leq C \inf_{\hat{s}} \sup_{s \in \mathcal{P}(R, \alpha, p)} \mathbb{E}_s \|s - \hat{s}\|^2.$$

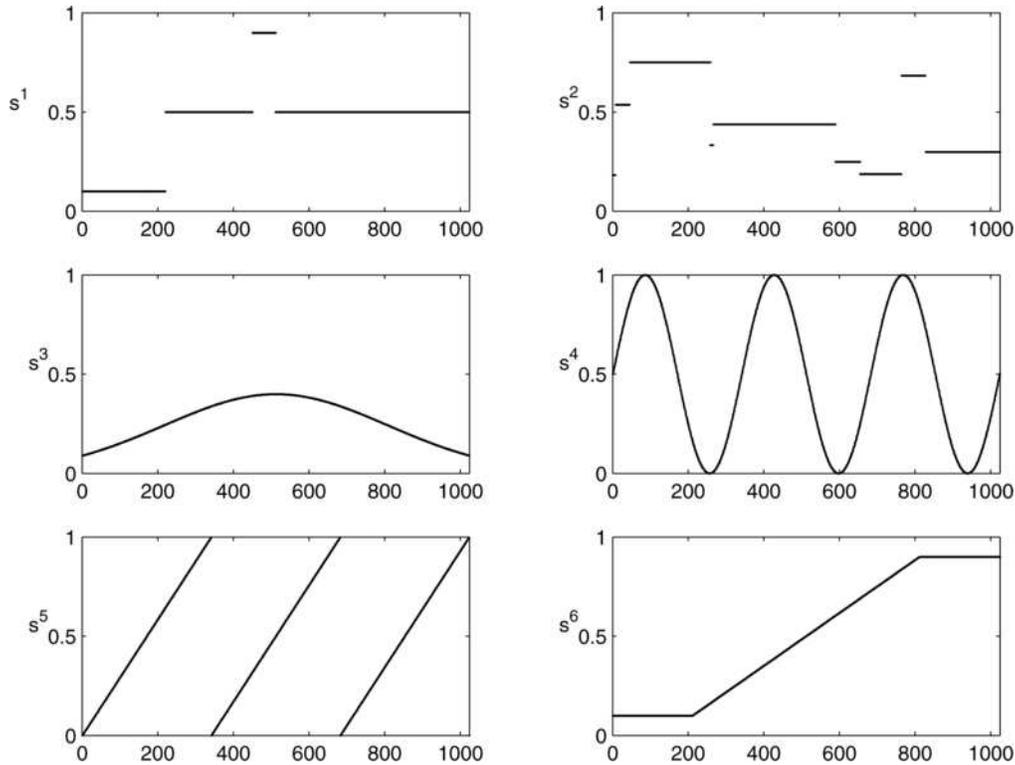

**Figure 1.** Points $(i, s_i^{(1)})$ of matrices $s^l, l = 1, \ldots, 6$.

## 6. Simulation study

In this section, numerical simulations are reported to compare the performances of the strategies EH, NEH and EI.

### 6.1. The simulation experiment

First, we fix $n = 2^{10} = 1024$, $r = 2$ and we consider six different matrices $s$ denoted by $s^l$, $l = 1, \ldots, 6$ and plotted in Figure 1. Precisely, only the points $(i, s_i^{(1)})$, $i = 1, \ldots, n$, corresponding to the first line of a given $s$ are plotted on this figure, since the second line is given by $s_i^{(2)} = 1 - s_i^{(1)}$. The considered functions $s$ are either piecewise constant (for $l = 1, 2$), or regular (for $l = 3, 4$), or piecewise linear (for $l = 5, 6$). Second, we fix $n = 2^{10}$, $r = 4$ and we consider two matrices $s$ that are constructed from the six previous ones. They are denoted, respectively, $s^7$ and $s^8$ and their four lines are plotted in Figure 2. For every considered $(n, r, s)$, we generate a sequence $(Y_1, \ldots, Y_n)$ and we compute the estimator $\tilde{s}$ with each of the three strategies. The obtained estimators are plotted together with $s$ on Figure 3 for $s = s^2, s^3, s^5$ (here again, as $r = 2$, only the first line of $s$ and $\tilde{s}$ are represented).

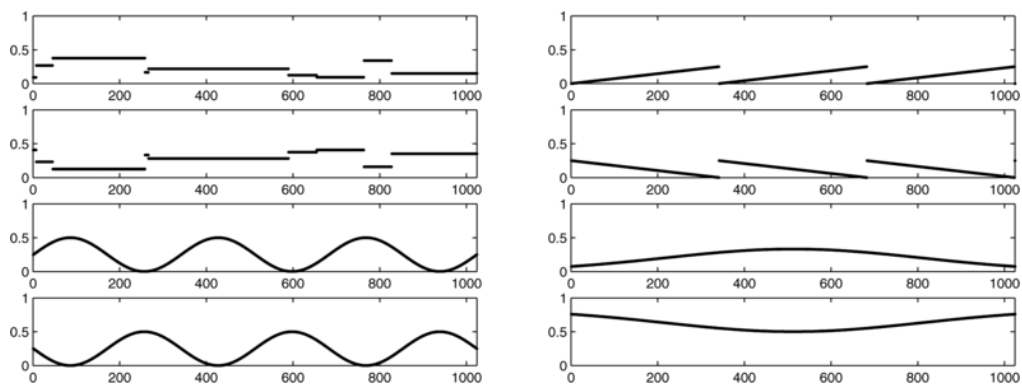

**Figure 2.** Points $(i, s_i^{(j)})$ of matrices $s^7$ and $s^8$.

We estimate the risk $\mathbb{E}_s[\|s - \tilde{s}\|^2]$ for each $s$ and each strategy using the following iterative method. At step $h$, we generate a sequence $(Y_1^h, \ldots, Y_n^h)$ according to $s$, we compute the estimator $\tilde{s}^h$ from this sequence and compute the average of the values $\|s - \tilde{s}^k\|^2$, $k = 1, \ldots, h$; we increment $h$ and repeat this step until the difference between two successive averages is lower than $10^{-2}$. Let $h^*$ be the number of performed steps. Then the risk is estimated by the empirical risk

$$\frac{1}{h^\star} \sum_{k=1}^{h^\star} \|s - \tilde{s}^k\|^2.$$

The considered strategies are based on a penalty function that calls in either two constants $c_1, c_2$ or a constant $c$, see (9) and (14). These constants have to be properly chosen so that the method performs well. So in order to compare the strategies in an objective way, we vary the constant(s) on a grid, estimate the estimator risk in each case and retain the optimal value $(c_1^*, c_2^*)$ or $c^*$ that leads to the smallest risk, which we call the minimal risk. Precisely, we vary $c_1$ from 0 to 1 by steps 0.1, $c_2$ from 0 to 6 by steps 0.1 and $c$ from 0 to 4 by steps 0.1. Table 1 gives the estimated minimal risk for each strategy and each considered $s$. Moreover, in order to study the influence of the constant(s), we plotted the risks obtained with different values of $(c_1, c_2)$ or $c$. As an example the graphs obtained for the estimation of $s^1$ using the EH and NEH strategies are presented on Figure 4.

To conclude the calibration study, let us notice that we can (at least heuristically) transpose to the NEH strategy the calibration method proposed in [6], since the NEH strategy calls in only one constant $c$. By doing so, our aim is to provide a data-driven calibration method taking into account the fact that the optimal choice for $c$ depends on the unknown function $s$. We use here the notation of Section 5. Given a value of the penalty constant $c$, the selected $J$, denoted here by $\hat{J}_c$, is computed by running the NEH algorithm with $0 \leq J \leq J_{\max}$ where $J_{\max}$ is fixed here to 7. We repeat the procedure for different values of $c$ increasing from 0, until $\hat{J}_c$ equals 0 (see [14] for similar computing). Then we consider the particular value of $c$, denoted $\hat{c}$, for which the difference between the dimensions of two successive selected models is the biggest. Finally, the retained penalty constant is $2\hat{c}$. Table 2 gives the mean and the standard deviation of the retained value for $c$, denoted by $\tilde{c}$ and $\tilde{\sigma}$, respectively, together with the estimated risk. For ease of comparison, we reported also in this table the minimal risk and its corresponding optimal constant $c^*$ defined above.

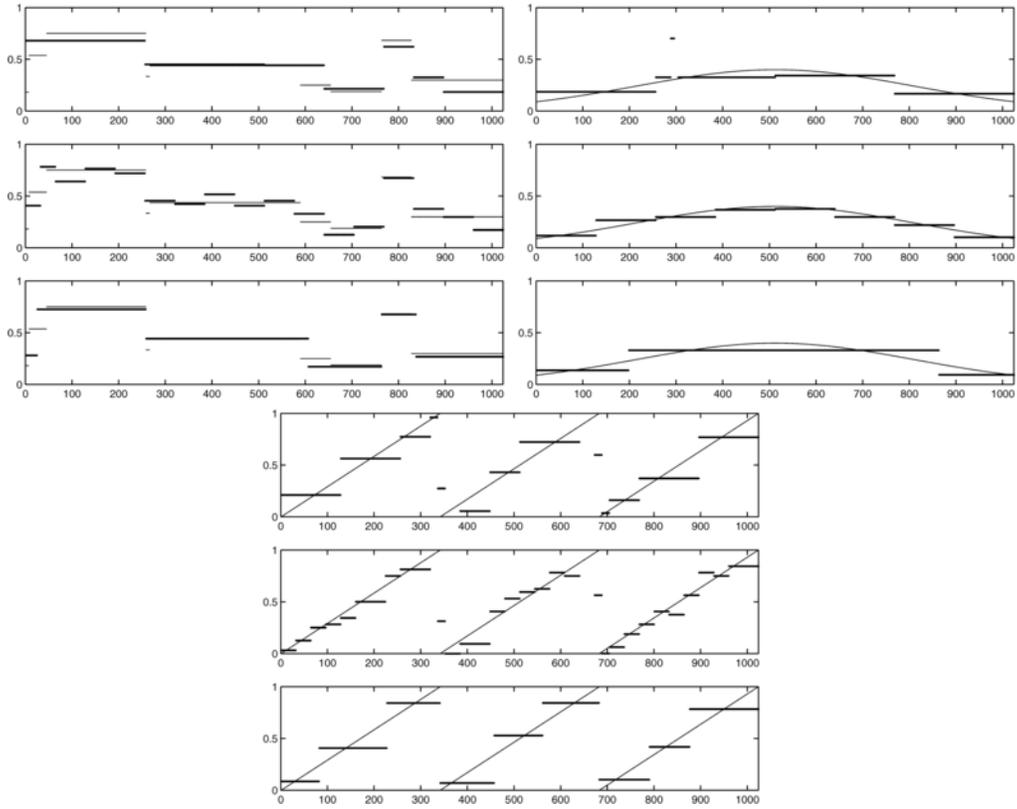

**Figure 3.** Points $(i, s_i^{(1)})$ (fine point) for $i = 1, \ldots, n$ and their estimators (thick point) for $s = s^2, s^3, s^5$ obtained, respectively, with the EH (top), NEH (middle) and EI (bottom) strategies.

**Table 1.** Minimal empirical risks for each strategy

|       | EH    | NEH   | EI    |
|-------|-------|-------|-------|
| $s^1$ | 14.48 | 10.50 |  5.34 |
| $s^2$ | 23.35 | 16.48 | 12.98 |
| $s^3$ |  7.52 |  5.52 |  8.03 |
| $s^4$ | 32.81 | 17.16 | 35.07 |
| $s^5$ | 33.19 | 28.10 | 27.46 |
| $s^6$ |  9.04 |  5.82 | 10.89 |
| $s^7$ | 27.84 | 21.17 | 29.84 |
| $s^8$ | 24.98 | 19.91 | 23.26 |

**Table 2.** Empirical risk $Risk_{\tilde{c}}$ obtained with the data-driven calibration and minimal empirical risks $Risk_{c^*}$ for the NEH strategy

| s | $\tilde{c}$ | $\tilde{\sigma}$ | $Risk_{\tilde{c}}$ | $c^*$ | $Risk_{c^*}$ |
|---|---|---|---|---|---|
| $s^1$ | 1.06 | 0.12 | 11.01 | 1.1 | 10.50 |
| $s^2$ | 1.06 | 0.10 | 16.52 | 1.1 | 16.48 |
| $s^3$ | 0.90 | 0.10 | 5.42 | 1.1 | 5.52 |
| $s^4$ | 0.74 | 0.10 | 16.99 | 0.8 | 17.16 |
| $s^5$ | 0.94 | 0.16 | 29.22 | 1.1 | 28.10 |
| $s^6$ | 0.82 | 0.10 | 6.71 | 1 | 5.82 |
| $s^7$ | 0.76 | 0.06 | 20.83 | 1.8 | 21.17 |
| $s^8$ | 0.70 | 0.05 | 20.17 | 1.6 | 19.91 |

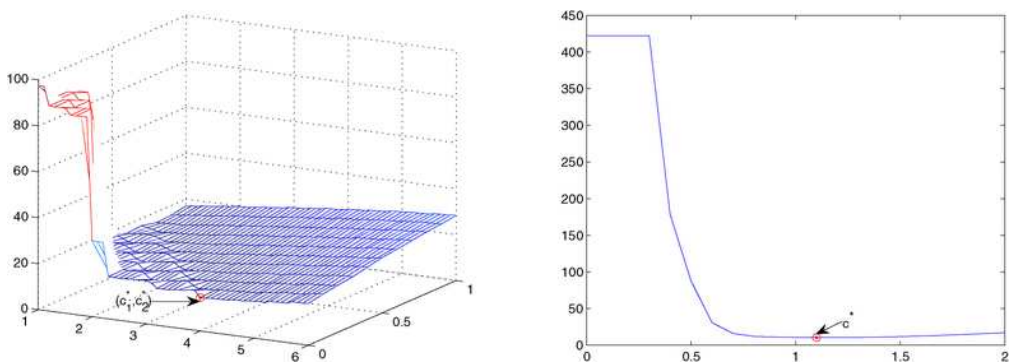

**Figure 4.** Estimated risks according to the different values of the penalty constants for the EH and NEH strategies, respectively, when estimating $s^1$.

### 6.2. Comments

In terms of risk, the NEH strategy outperforms the EH strategy for every considered $s$, see Table 1. This is consistent with our theoretical results since the NEH strategy satisfies oracle and minimax results where the EH strategy satisfies these results only up to a log factor. The EI strategy outperforms the NEH strategy, as expected, when $s$ is piecewise constant, whereas the NEH strategy outperforms the EI strategy for all the other considered $s$.

In terms of computational issues, both strategies based on Haar functions outperform the EI strategy. Indeed, with these strategies the computational effort is only of order $\mathcal{O}(n \log n)$, see Sections 4 and 5. On the contrary, even using dynamic programming the computational complexity of the EI strategy is of order $\mathcal{O}(n^3)$.

The last issue concerns the choice of the penalty constants. First, we observed in the simulation results that for every considered $s^l$ and whatever the strategy, small departures from the optimal values of the penalty constants do not damage the risk of the estimator (see Figure 4 for the estimation of $s^1$ with the EH and NEH strategies). Now a great advantage of the NEH strategy in comparison with both exhaustive strategies comes from the fact that it calls in only one constant. One can then use the data-driven calibration method described in Section 6.1 and we observe that this method works successfully since the obtained risks are very close to the minimal risks (see Table 2).

## 7. Application to change point detection

In this section, we provide an application of our method to the change point detection problem. Precisely, we adapt our method to this problem by providing an algorithm that combines the NEH and the EI strategies, then we run this algorithm on a DNA sequence.

### 7.1. An hybrid algorithm

Our aim is to build a piecewise constant estimator of $s$ in such a way that the change points in the estimator reflect change points in the underlying distribution $s$, when the sample is large. For this task, we first estimate $s$ using the NEH strategy and the penalty function

$$\text{pen}(m) = 2\hat{c}D_m,$$

where $\hat{c}$ is obtained with the data-driven calibration method described in Section 6.1. The obtained estimator $\tilde{s}$ is piecewise constant and properly estimates $s$ but, due to its construction, it possesses many jumps and not all of them have to be interpreted as change points in $s$. Thus, we use the EI strategy to decide which change points in $\tilde{s}$ should be interpreted as change points in $s$. To be more precise, let us denote by $\hat{\mathcal{J}}$ the set of points $i \in \{2, \ldots, n\}$ such that $\tilde{s}_i \neq \tilde{s}_{i-1}$. We perform the EI strategy on the collection of all the partitions of $\{1,\ldots,n\}$ built on $\hat{\mathcal{J}}$. Here, the penalty function takes the form

$$\mathrm{pen}_1(m) = D_m\left(c_1 \log\left(\frac{|\hat{\mathcal{J}}|}{D_m}\right) + c_2\right), \tag{15}$$

where $c_1$ and $c_2$ have to be chosen and where $|\hat{\mathcal{J}}|$ denotes the cardinality of $\hat{\mathcal{J}}$. This penalty calls in two constants that may pose difficulties on practical applications, so we slightly modify the EI strategy by considering a simpler form of penalty function. It is easy to see from Section 3 that a penalty of the form $cD_m$ is allowed in connection with the collection of models considered there, provided $c \geq k\log(n)$ with some large enough $k$. Thus, instead of (15), we consider the penalty function

$$\mathrm{pen}_2(m) = cD_m. \tag{16}$$

**Table 3.** Minimal empirical risks for the hybrid and the EI strategies

|       | NEH–EI$_1$ | NEH–EI$_2$ | EI    |
|-------|------------|------------|-------|
| $s^1$ | 5.34       | 5.06       | 5.34  |
| $s^2$ | 14.5       | 15.48      | 12.98 |

To summarize, our hybrid algorithm consists in two steps: Provide with the NEH strategy a piecewise constant estimator with many jumps, then apply the EI strategy. This is of interest when the sample is large, since the EI strategy (and the segmentation methods described in Section 1) has a computational complexity that is too high to be implemented on a large sequence of length $n$, but can typically be implemented if the set of possible change points is restricted to $\hat{\mathcal{J}}$. Such an algorithm is similar in spirit to the algorithms in Gey and Lebarbier [11] and Akakpo [2]. In [11], the first step estimator is built using the CART (Classification and Regression Tree) algorithm but suffers from a higher computational complexity. In [2], it is inspired by our own method and relies on dyadic intervals.

We report below a small simulation study to justify the use of (16). We denote by NEH–EI$_1$ and NEH–EI$_2$ the hybrid strategies with a penalty of the form (15) and (16), respectively, at the EI step. We vary $c_1$, $c_2$ and $c$ on the same grids as in Section 6.1 and compute the minimal empirical risks obtained for both hybrid strategies NEH–EI$_1$ and NEH–EI$_2$, when $s$ is either $s^1$ or $s^2$. The results are given in Table 3, where we report also the minimal empirical risks obtained in Section 6.1 using the EI strategy (recall that among the strategies studied there, EI provides the smallest risk on $s^1$ and $s^2$). The risks obtained with the three strategies are comparable. This justifies the use of the hybrid strategies: They can be implemented even on large samples (whereas EI cannot) and compares, in terms of risk, to EI on the considered examples.

### 7.2. Application on a DNA sequence

In this section, we run the NEH–EI$_2$ strategy on a DNA sequence and compare our estimated change points with the annotation contained in the GenBank database. Precisely, we consider the *Bacillus subtilis* sequence reduced to $n = 2^{21} = 2{,}097{,}152$ bases in length. We consider this DNA sequence as a sequence of independent multinomial variables with values in the alphabet $\{A, C, G, T\}$ of length $r = 4$, and we run NEH–EI$_2$ with penalty constants calibrated with the data-driven method of Section 6.1. This provides us an estimator of $s$. We plot the four lines of this estimator together with the GenBank annotation using the software MuGeN [12]. For the sake of readability, Figure 5 shows the results for the first 178,000 bases only. The arrows correspond to genes and their direction to the transcription direction. All the regions composed of genes coding for ribosomal RNA (dark genes on Figure 5) are detected by NEH–EI$_2$. Moreover NEH–EI$_2$ detects some changes of transcription direction. These results are close to those obtained in [13] on the first 200,000 bases of the sequence. However our strategy is much faster and allows the study of a much greater part of the sequence.

## 8. Proof of the theorems

Throughout the section, $C$ denotes a constant that may change from line to line. The cardinality of a set $E$ is denoted by $|E|$. Moreover, we set for every $x \geq 0$,

$$\lfloor x \rfloor = \sup\{l \in \mathbb{N}, l \leq x\} \quad \text{and} \quad \lceil x \rceil = \inf\{l \in \mathbb{N}, l \geq x\}.$$

### 8.1. Proof of Theorem 1

In order to prove Theorem 1, we follow the line of the proof of Theorem 2 in [5]. Fix $m \in \mathcal{M}$. By definition of $\tilde{s}$ and $\hat{s}_m$ we have

$$\gamma_n(\tilde{s}) + \text{pen}(\hat{m}) \leq \gamma_n(\bar{s}_m) + \text{pen}(m).$$

But

$$\gamma_n(t) = -2\langle \varepsilon, t \rangle + \|s - t\|^2 - \|s\|^2$$

for all $r \times n$ matrices $t$. Therefore

$$\|s - \tilde{s}\|^2 \leq \|s - \bar{s}_m\|^2 + 2\langle \varepsilon, \tilde{s} - \bar{s}_m \rangle + \text{pen}(m) - \text{pen}(\hat{m}). \tag{17}$$

We need a concentration inequality to control the random term $\langle \varepsilon, \tilde{s} - \bar{s}_m \rangle$. More precisely, we will control $\langle \varepsilon, \hat{s}_{m'} - \bar{s}_m \rangle$ uniformly in $m' \in \mathcal{M}$ by using the following lemma (see Appendix A.1 for a proof).

Finite distribution estimation via model selection 495

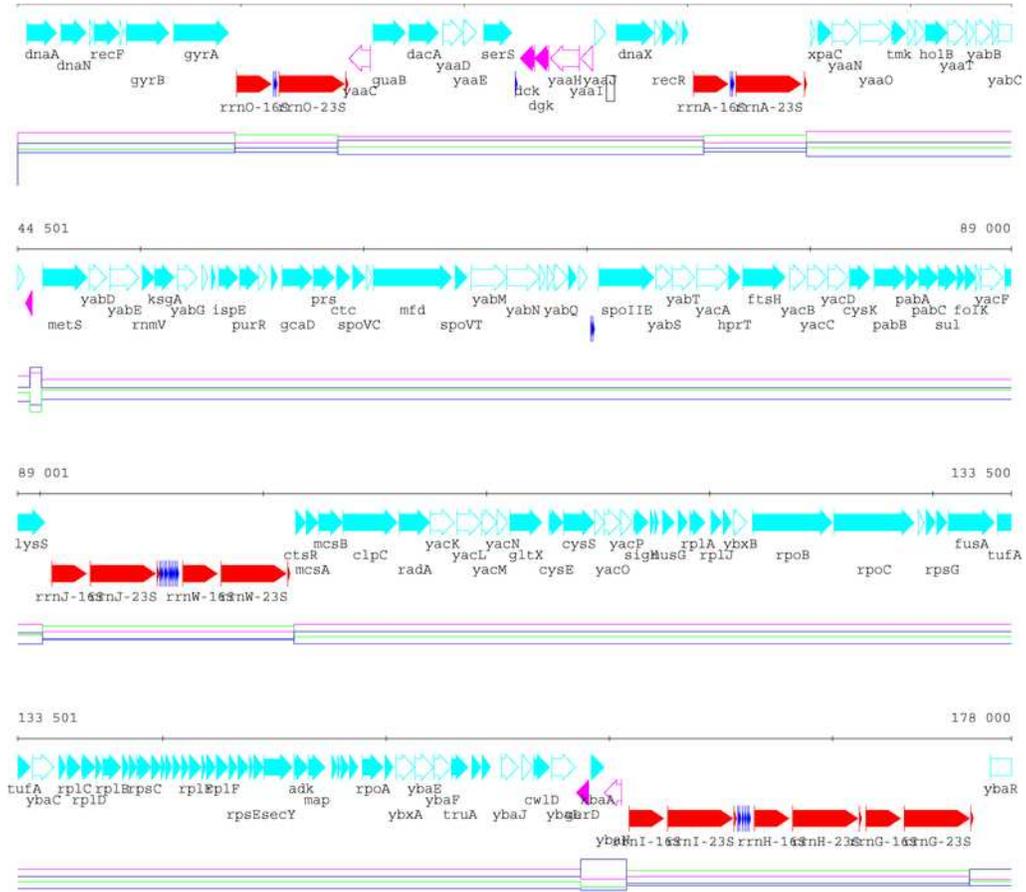

**Figure 5.** Annotation of the first 178 000 bases of *Bacillus subtilis*, and estimator obtained with NEH-EI$_2$.

**Lemma 1.** *For every $m' \in \mathcal{M}$ and $x \geq 0$, we have*

$$\mathbb{P}_s\left(\sup_{u \in \mathbb{R}^r \otimes (S_m + S_{m'})} \frac{\langle \varepsilon, u \rangle}{\|u\|} \geq \sqrt{D_m + D_{m'}} + 2\sqrt{2x}\right) \leq \exp(-x).$$

Fix $\xi > 0$ and let $\Omega_\xi(m)$ denote the event

$$\Omega_\xi(m) = \bigcap_{m' \in \mathcal{M}} \left\{\sup_{u \in \mathbb{R}^r \otimes (S_m + S_{m'})} \frac{\langle \varepsilon, u \rangle}{\|u\|} \leq \sqrt{D_m + D_{m'}} + 2\sqrt{2(L_{m'} D_{m'} + \xi)}\right\}.$$

It follows from the definition of $\Sigma$ that

$$\mathbb{P}_s(\Omega_\xi(m)) \geq 1 - \Sigma \exp(-\xi). \tag{18}$$

Moreover, (17) proves that on $\Omega_\xi(m)$,

$$\|s-\tilde{s}\|^2 \leq \|s-\bar{s}_m\|^2 + 2(\sqrt{D_m + D_{\hat{m}}} + 2\sqrt{2(L_{\hat{m}}D_{\hat{m}} + \xi)}) \times \|\tilde{s} - \bar{s}_m\| + \text{pen}(m) - \text{pen}(\hat{m}).$$

Note that $\|\tilde{s} - \bar{s}_m\| \leq \|\tilde{s} - s\| + \|s - \bar{s}_m\|$ and that for every positive number $a$, $b$ and $\theta$ one has $2\sqrt{ab} \leq \theta a + \theta^{-1} b$. Applying this inequality twice, we get that on $\Omega_\xi(m)$ and for any $\eta \in (0,1)$,

$$\|s-\tilde{s}\|^2 \leq \|s-\bar{s}_m\|^2 + (1-\eta)((1+\eta)\|\tilde{s}-s\|^2 + (1+\eta^{-1})\|s-\bar{s}_m\|^2)$$
$$+ \frac{1}{1-\eta}(\sqrt{D_m + D_{\hat{m}}} + 2\sqrt{2(L_{\hat{m}}D_{\hat{m}} + \xi)})^2 + \text{pen}(m) - \text{pen}(\hat{m}).$$

But for all positive numbers $a$ and $b$, one has $\sqrt{a+b} \leq \sqrt{a} + \sqrt{b}$; so for any $\eta \in (0,1)$ we get

$$(\sqrt{D_m + D_{\hat{m}}} + 2\sqrt{2(L_{\hat{m}}D_{\hat{m}} + \xi)})^2 \leq (\sqrt{D_m} + \sqrt{D_{\hat{m}}}(1 + 2\sqrt{2L_{\hat{m}}}) + 2\sqrt{2\xi})^2$$
$$\leq 2(1+\eta^{-1})(D_m + 8\xi) + (1+\eta)D_{\hat{m}}(1 + 2\sqrt{2L_{\hat{m}}})^2.$$

In particular, for $\eta$ such that

$$K = \frac{1+\eta}{1-\eta},$$

that is, $\eta = (K-1)/(K+1)$, it follows from the assumption on the penalty function pen that on $\Omega_\xi(m)$,

$$\left(\frac{K-1}{K+1}\right)^2 \|s-\tilde{s}\|^2 \leq \frac{K^2 + 4K - 1}{K^2 - 1}\|s-\bar{s}_m\|^2 + \frac{2K(K+1)}{K-1}(D_m + 8\xi) + \text{pen}(m).$$

But $K^2 + 4K - 1 \leq 2K(K+1)$ and $D_m \leq \text{pen}(m)/K$ so on $\Omega_\xi(m)$,

$$\|s-\tilde{s}\|^2 \leq \frac{2K(K+1)^2}{(K-1)^3}\{\|s-\bar{s}_m\|^2 + \text{pen}(m)\} + 16\frac{K(K+1)^3}{(K-1)^3}\xi.$$

It then follows from the definition of $\Omega_\xi(m)$ and (18) that there exists a random variable $V \geq 0$ such that on the one hand,

$$\|s-\tilde{s}\|^2 \leq \frac{2K(K+1)^2}{(K-1)^3}\{\|s-\bar{s}_m\|^2 + \text{pen}(m)\} + 16\frac{K(K+1)^3}{(K-1)^3}V$$

and on the other hand,

$$\mathbb{P}_s(V > \xi) \leq \Sigma \exp(-\xi) \qquad \forall \xi > 0.$$

Integrating the latter inequality yields $\mathbb{E}_s(V) \leq \Sigma$, which completes the proof of the theorem since $m$ is arbitrary.

## 8.2. Proof of Theorems 2 and 3

We need more notation. Let $P_s$ be the distribution of $X$ under $\mathbb{P}_s$. Let $u$ be the $r \times n$ matrix such that $u^{(j)}$ is the vector of $\mathbb{R}^n$ with all components equal to $1/2$ for $j \in \{1, 2\}$, and in the case where $r > 2$, $u^{(j)}$ is the null vector of $\mathbb{R}^n$ for all $j > 2$. If there exists $N \in \mathbb{N}$ such that $n = 2^N$, we denote by $\{\phi_\lambda, \lambda = 1, \ldots, n\}$ the Haar basis given in Definition 1 with indexes arranged in lexical order: an element $\lambda \in \Lambda(j)$ is written $\lambda = 1$ if $j = -1$ and $\lambda = 1 + l$ for some $l \in \{2^j, \ldots, 2^{j+1} - 1\}$ if $j = 0, \ldots, N - 1$. In that case, for every $D_0 \in \{1, \ldots, n-1\}$, $D \in \{1, \ldots, D_0\}$ and $A > 0$, we denote by $\tilde{\mathcal{P}}(A, D_0, D)$ the set of those $r \times n$ matrices $s$ such that each column of $s$ is a probability on $\{1, \ldots, r\}$ and $s = u + t$ where for all $j$

$$t^{(j)} = \sum_{\lambda \in m} \beta_\lambda^{(j)} \phi_\lambda,$$

$m$ is a subset of $\{2, \ldots, 1 + D_0\}$ with cardinality $D$ and $\|\beta_\lambda\|_r \leq A$ for all $\lambda \in m$. We will see later that when $n = 2^N$ for some $N \in \mathbb{N}$,

$$\tilde{\mathcal{P}}(A, D_0, D) \subset \mathcal{P},$$

where $\mathcal{P}$ denotes either $\mathcal{P}_k$ or $\mathcal{P}(R, \alpha, p)$, $A$ is a well-chosen positive number and $D_0 \geq D$ are well-chosen integers. It follows that when $n = 2^N$,

$$\sup_{s \in \mathcal{P}} \mathbb{E}_s \|s - \hat{s}\|^2 \geq \sup_{s \in \tilde{\mathcal{P}}(A, D_0, D)} \mathbb{E}_s \|s - \hat{s}\|^2.$$

In order to get a lower bound on $\mathcal{P}$, we thus first compute a lower bound on $\tilde{\mathcal{P}}(A, D_0, D)$ for arbitrary $A$, $D_0$ and $D$ (see Appendix A.2 for a proof).

**Lemma 2.** *Assume $n = 2^N$. There exists a universal constant $C > 0$ such that*

$$\inf_{\hat{s}} \sup_{s \in \tilde{\mathcal{P}}(A, D_0, D)} \mathbb{E}_s \|s - \hat{s}\|^2 \geq CD \left( A^2 \wedge \left( \log\left(\frac{D_0}{D}\right) + 1 \right) \wedge \frac{n}{D_0} \right) \tag{19}$$

*for every $A$, $D$ and $D_0$ with $A > 0$ and $1 \leq D \leq D_0 < n$.*

We will also use the following two lemmas.

**Lemma 3.** *Let $u$ and $v$ be $r \times n$ matrices the columns of which are probabilities on $\{1, \ldots, r\}$. Then,*

$$\inf_{\hat{s}} \sup_{s \in \{u, v\}} \mathbb{E}_s \|s - \hat{s}\|^2 \geq \frac{\|u - v\|^2}{4} \left( 1 - \sqrt{\frac{K(P_u, P_v)}{2}} \right),$$

where $K(\cdot,\cdot)$ denotes Kullback–Leibler divergence.

**Proof.** For every estimator $\hat{s}$,

$$\sup_{s\in\{u,v\}} \mathbb{E}_s\|s-\hat{s}\|^2 \geq \frac{1}{2}\int(\|u-\hat{s}\|^2 + \|v-\hat{s}\|^2)\min(\mathrm{d}P_u,\mathrm{d}P_v).$$

By convexity of the function $x \mapsto x^2$ combined to the triangle inequality,

$$\|u-\hat{s}\|^2 + \|v-\hat{s}\|^2 \geq \|u-v\|^2/2.$$

Moreover, Pinsker's inequality yields

$$\int \min(\mathrm{d}P_u,\mathrm{d}P_v) \geq 1 - \sqrt{\frac{K(P_u,P_v)}{2}},$$

so we get the result. □

**Lemma 4.** *Let $u$ and $v$ be $r \times n$ matrices such that each column is a probability on $\{1,\ldots,r\}$ and each component either equals zero or belongs to $[1/4,1]$. If $P_u \ll P_v$ and $P_v \ll P_u$, then*

$$K(P_u,P_v) \leq 4\|u-v\|^2,$$

*where $K(\cdot,\cdot)$ denotes Kullback–Leibler divergence.*

**Proof.** The result is obtained by bounding the Kullback–Leibler divergence by the $\chi^2$ divergence:

$$K(P_u,P_v) \leq \sum_{i=1}^n \sum_{j=1}^r v_i^{(j)}\left(\frac{u_i^{(j)}}{v_i^{(j)}} - 1\right)^2 \leq 4\|u-v\|^2. \qquad \Box$$

**Proof of Theorem 2.** We distinguish between three cases. First, we consider the case where $k \geq 4$, since in that case we can apply Lemma 2. Then we assume $k = 3$. If $n$ is bounded ($n < 12$, say) then the log term in the lower bound can be viewed as a constant and we will see that the result is an easy consequence of Lemma 3. In the last case ($k = 3$ and $n$ possibly large), we use another argument to get the log term in the lower bound.

Case 1. $n \geq k \geq 4$. Note first that it suffices to prove the result under the additional assumption that there exists an integer $N$ such that $n = 2^N$. Indeed, assume there exists an absolute $C > 0$ such that

$$\inf_{\hat{s}} \sup_{s \in \mathcal{P}_k} \mathbb{E}_s\|s-\hat{s}\|^2 \geq Ck\left(\log\left(\frac{n}{k}\right) + 1\right)$$

for all $k \in \{4,\ldots,n\}$, provided $n$ is a power of 2. Then, in the general case, there exists an integer $N \geq 2$ such that $2^N \leq n \leq 2^{N+1}$ and it is easy to derive from the last display

that

$$\inf_{\hat{s}} \sup_{s \in \mathcal{P}_k} \mathbb{E}_s \|s - \hat{s}\|^2 \geq Ck\left(\log\left(\frac{2^N}{k}\right) + 1\right)$$

for every $k \in \{4, \ldots, 2^N\}$. But $2^N \geq n/2$ so (possibly reducing the constant $C$) we obtain the result for every $k \in \{4, \ldots, 2^N\}$. Now, every $k$ with $2^N < k \leq n$ satisfies $k \leq 2k'$, where $k' = 2^N$, and we have $\mathcal{P}_{k'} \subset \mathcal{P}_k$. Hence

$$\inf_{\hat{s}} \sup_{s \in \mathcal{P}_k} \mathbb{E}_s \|s - \hat{s}\|^2 \geq \inf_{\hat{s}} \sup_{s \in \mathcal{P}_{k'}} \mathbb{E}_s \|s - \hat{s}\|^2 \geq C\frac{k}{2}\left[\log\left(\frac{n}{k}\right) + 1\right]$$

and we get the required result for every $n \geq 4$ and $k \in \{4, \ldots, n\}$. We thus assume in the sequel that $n = 2^N$ and we fix $k \in \{4, \ldots, n\}$. We then have $\tilde{\mathcal{P}}(A, D_0, D) \subset \mathcal{P}_k$, provided $k \geq 4D$. Using Lemma 2 with $D = \lfloor k/4 \rfloor$ and $A = \sqrt{n/D_0}$ thus yields

$$\inf_{\hat{s}} \sup_{s \in \mathcal{P}_k} \mathbb{E}_s \|s - \hat{s}\|^2 \geq C\lfloor k/4 \rfloor \left(\left(\log\left(\frac{D_0}{\lfloor k/4 \rfloor}\right) + 1\right) \wedge \frac{n}{D_0}\right)$$

for every $D_0 \geq \lfloor k/4 \rfloor$ with $D_0 < n$. But $\lfloor k/4 \rfloor \geq k/7$ and $\lfloor k/4 \rfloor \leq k/4$, so

$$\inf_{\hat{s}} \sup_{s \in \mathcal{P}_k} \mathbb{E}_s \|s - \hat{s}\|^2 \geq \frac{Ck}{7}\left(\left(\log\left(\frac{4D_0}{k}\right) + 1\right) \wedge \frac{n}{D_0}\right). \tag{20}$$

Let $D_0$ be an integer that makes the right-hand side of the previous inequality as large as possible, that is, some integer that approximately equates the terms $\log(4D_0/k)$ and $n/D_0$:

$$D_0 = \left\lfloor \frac{n}{\log(4n/k)} \right\rfloor.$$

We have $-x \log x \leq e^{-1}$ for every $x > 0$, so

$$\frac{k}{4n}\log\left(\frac{4n}{k}\right) \leq \frac{1}{2}, \tag{21}$$

which implies $D_0 \geq \lfloor k/4 \rfloor$. Moreover, $D_0 < n$ since $k \leq n$, so (20) holds for this value of $D_0$. By (21) and the definition of $D_0$,

$$\frac{n}{2\log(4n/k)} \leq D_0 \leq \frac{n}{\log(4n/k)},$$

so (20) proves that there exists an absolute constant $C > 0$ such that

$$\inf_{\hat{s}} \sup_{s \in \mathcal{P}_k} \mathbb{E}_s \|s - \hat{s}\|^2 \geq \frac{Ck}{7}\left(\left(\log\left(\frac{2n}{k}\right) + 1 - \log\log\left(\frac{4n}{k}\right)\right) \wedge \left(\log\left(\frac{n}{k}\right) + \log 4\right)\right).$$

Straightforward computations then yield the result.

Case 2. $k = 3$ and $n \geq 12$. Let

$$L = \left\lfloor \frac{n}{\lfloor \log n \rfloor} \right\rfloor$$

and for every $l = 0, \ldots, L-1$, let $\varphi_l$ be that vector of $\mathbb{R}^n$ with the $i$th component $\varphi_{li} = 1/\sqrt{\lfloor \log n \rfloor}$ if $i \in \{1 + l\lfloor \log n \rfloor, \ldots, (l+1)\lfloor \log n \rfloor\}$ and $\varphi_{li} = 0$ otherwise. For every $l = 0, \ldots, L-1$, let $s_l$ be the $r \times n$ matrix defined by

$$\begin{cases} s_l^{(1)} = 1/2 + a\varphi_l, \\ s_l^{(2)} = 1/2 - a\varphi_l, \\ s_l^{(j)} = 0, & \text{for } j > 2 \text{ (in the case } r > 2\text{)}, \end{cases}$$

where

$$a = \sqrt{\frac{\lfloor \log n \rfloor}{48}}.$$

The system $\{\varphi_l, l = 0, \ldots, L-1\}$ is orthonormal so $\|s_l - s_{l'}\|^2 = 4a^2$ for every $l \neq l'$. Moreover, $s_{li}^{(j)} \geq 1/4$ for every $i = 1, \ldots, n$ and $j = 1, 2$ so Lemma 4 yields

$$K(P_{s_l}, P_{s_{l'}}) \leq 4\|s_l - s_{l'}\|^2 \leq 16a^2$$

for every $l \neq l'$. We have $L \geq 6$ and $16a^2/\log L \leq 2/3$, since $n \geq 12$ and $L \geq \sqrt{n}$. By Proposition 9 of [5] with $\mathcal{C} = \{s_l, l = 0, \ldots, L-1\}$, $\eta = 2a$ and $H = 16a^2$,

$$\inf_{\hat{s}} \sup_{s \in \mathcal{C}} \mathbb{E}_s \|s - \hat{s}\|^2 \geq \frac{1}{3 \times 48} \lfloor \log n \rfloor,$$

where the infimum is taken over the estimators $\hat{s}$ such that each column $\hat{s}^{(j)}$ is a probability on $\{1, \ldots, r\}$. Since $\mathcal{C} \subset \mathcal{P}_3$, we get the result by extending the infimum over all possible estimators.

Case 3. $k = 3$ and $n < 12$. Let $v \in \mathbb{R}^r \otimes \mathbb{R}^n$ be defined by $v_1^{(1)} = 3/4$, $v_1^{(2)} = 1/4$ and $v_i^{(j)} = u_i^{(j)}$ for every $(i, j) \neq (1, 1), (1, 2)$, where $u$ is defined at the beginning of Section 8.2. Then $\|u - v\|^2 = 1/8$ and Lemma 4 yields $K(P_u, P_v) \leq 1/2$. It thus follows from Lemma 3 that

$$\inf_{\hat{s}} \sup_{s \in \{u,v\}} \mathbb{E}_s \|s - \hat{s}\|^2 \geq \frac{1}{64}.$$

Both $u$ and $v$ belong to $\mathcal{P}_3$, so there exists an absolute constant $C > 0$ such that

$$\inf_{\hat{s}} \sup_{s \in \mathcal{P}_3} \mathbb{E}_s \|s - \hat{s}\|^2 \geq C.$$

The result follows since $n < 12$ and $k = 3$. □

**Proof of Theorem 3.** For every $D \in \{1, \ldots, n-1\}$, it is easy to see that

$$\tilde{\mathcal{P}}(\sqrt{n}RD^{-(\alpha+1/2)}, D, D) \subset \mathcal{P}(R, \alpha, p).$$

Therefore,

$$\sup_{s \in \mathcal{P}(R,\alpha,p)} \mathbb{E}_s \|s - \hat{s}\|^2 \geq \sup_{s \in \tilde{\mathcal{P}}(\sqrt{n}RD^{-(\alpha+1/2)}, D, D)} \mathbb{E}_s \|s - \hat{s}\|^2$$

for every estimator $\hat{s}$ and $D \in \{1 \ldots, n-1\}$. It then follows from Lemma 2 that there exists an absolute constant $C > 0$ such that for every $D \in \{1, \ldots, n-1\}$,

$$\inf_{\hat{s}} \sup_{s \in \mathcal{P}(R,\alpha,p)} \mathbb{E}_s \|s - \hat{s}\|^2 \geq CD((nR^2 D^{-(2\alpha+1)}) \wedge 1).$$

We choose $D$ that approximately equates the term $nR^2 D^{-(2\alpha+1)}$ with 1; we set

$$D = \lfloor (nR^2)^{1/(2\alpha+1)} \rfloor.$$

Then $D \geq 1$ since $nR^2 > 1$ and $D < n$ since $R < n^\alpha$, so we get

$$\inf_{\hat{s}} \sup_{s \in \mathcal{P}(R,\alpha,p)} \mathbb{E}_s \|s - \hat{s}\|^2 \geq C \lfloor (nR^2)^{1/(2\alpha+1)} \rfloor.$$

The result now follows from

$$\lfloor (nR^2)^{1/(2\alpha+1)} \rfloor \geq \tfrac{1}{2} (nR^2)^{1/(2\alpha+1)}. \qquad \square$$

### 8.3. Proof of Theorems 4 and 5

In order to prove Theorems 4 and 5, we first control the approximation terms

$$\inf_{m \in \mathcal{M}, D_m = D} \|s - \bar{s}_m\|^2$$

for $s \in \mathcal{P}(R, \alpha, p)$, where $\bar{s}_m$ is the orthogonal projection of $s$ onto $\mathbb{R}^r \otimes S_m$. For this task, we use the following lemma.

**Lemma 5.** *Let $R > 0$, $p \in (0, 2]$, $\alpha > 1/p - 1/2$ and $s \in \mathcal{P}(R, \alpha, p)$. Let $J \in \{0, \ldots, N-1\}$ and let $\mathcal{M}_J$ be defined as in Section 5. There exists some $C > 0$ that only depends on $p$ and $\alpha$ such that*

$$\inf_{m \in \mathcal{M}_J} \|s - \bar{s}_m\|^2 \leq CnR^2 2^{-2\alpha J}.$$

**Proof.** In the sequel, we use the notation of Definition 8 of [5], page 237. Moreover, we set

$$\beta_\lambda^{(j)} = \langle s^{(j)}, \phi_\lambda \rangle_n$$

for every $j \in \{1,\ldots,r\}$ and $\lambda \in \Lambda$. By definitions of $\bar{s}_m$ and $\mathcal{M}_J$, we have

$$\inf_{m \in \mathcal{M}_J} \|s - \bar{s}_m\|^2 = \sum_{k=0}^{N-J-1} \sum_{\lambda \notin \overline{\Lambda(J+k)}(|\Lambda'(J+k)|)} \|\beta_\lambda\|_r^2.$$

Since $s \in \mathcal{P}(R,\alpha,p)$, we have

$$\sum_{\lambda \in \Lambda(J+k)} \|\beta_\lambda\|_r^p \leq n^{p/2} R^p 2^{-p(J+k)(\alpha+1/2-1/p)}.$$

By Lemma 2 of [5], page 246 we thus have

$$\sum_{\lambda \notin \overline{\Lambda(J+k)}(|\Lambda'(J+k)|)} \|\beta_\lambda\|_r^2 \leq (|\Lambda'(J+k)| + 1)^{1-2/p} (\sqrt{n} R 2^{-(J+k)(\alpha+1/2-1/p)})^2.$$

But $1 - 2/p \leq 0$ and $|\Lambda'(J+k)| + 1 \geq 2^J (k+1)^{-3}$ so

$$\inf_{m \in \mathcal{M}_J} \|s - \bar{s}_m\|^2 \leq nR^2 2^{-2\alpha J} \sum_{k \geq 0} (k+1)^{-3(1-2/p)} 2^{-2k(\alpha+1/2-1/p)},$$

which proves the lemma. □

**Proof of Theorem 4.** Let $s$ be an arbitrary element of $\mathcal{P}(R,\alpha,p)$. All the elements $m$ of a given $\mathcal{M}_J$ have the same cardinality and satisfy $|m| \leq c_0 2^J$, with $c_0 = 1 + \sum_{k \geq 1} k^{-3}$. Let $D \in \{1,\ldots,n\}$ with $D \geq c_0$ and let

$$J = \sup\{j = 0,\ldots,N-1 \text{ s.t. } c_0 2^j \leq D\}.$$

We have $|m| \leq c_0 2^J \leq D$ for every $m \in \mathcal{M}_J$, so it follows from the definition of $\bar{s}_m$ that

$$\inf_{m \in \mathcal{M}, D_m = D} \|s - \bar{s}_m\|^2 \leq \inf_{m \in \mathcal{M}_J} \|s - \bar{s}_m\|^2.$$

By the definition of $J$ we have $D < c_0 2^{J+1}$, so it follows from Lemma 5 that

$$\inf_{m \in \mathcal{M}, D_m = D} \|s - \bar{s}_m\|^2 \leq (2c_0)^{2\alpha} C n R^2 D^{-2\alpha}.$$

If $c_1$ and $c_2$ are large enough, then (2) holds with $L_m = 2 + \log(n/D_m)$ and $\Sigma \leq 1$, see (7) and (8). It thus follows from Theorem 1 and the definition of the penalty function pen that there exists $C > 0$ that only depends on $\alpha$, $p$, $c_1$ and $c_2$ such that for every $D = 1,\ldots,n$ with $D \geq c_0$,

$$\mathbb{E}_s \|s - \tilde{s}\|^2 \leq C \left\{ nR^2 D^{-2\alpha} + D\left(1 + \log\left(\frac{n}{D}\right)\right) \right\}. \tag{22}$$

Possibly enlarging $C$, one gets (22) for every $D = 1, \ldots, n$. In order to optimize this inequality, we search for the $D$ that provides the smallest upper bound, that is, some $D$ that approximately equates the terms $nR^2 D^{-2\alpha}$ and $D \log(n/D)$. For this task let us define the function $f$ on $(1, n^{2\alpha+1})$ by

$$f(x) = x^{1/(2\alpha+1)} (\log(nx^{-1/(2\alpha+1)}))^{-1/(2\alpha+1)}$$

and let $D = \lceil f(nR^2) \rceil \geq 1$. The function $x \mapsto x \log x$ is increasing on $[e^{-1}, \infty)$ and, by assumption, $n^\alpha > R\sqrt{\rho}$. Therefore

$$\frac{n^{2\alpha}}{R^2} \log \frac{n^{2\alpha}}{R^2} > 2\alpha + 1,$$

which implies $D \leq n$. Moreover, $f$ is increasing on $(1, n^{2\alpha+1})$ and, by assumption, $1 < \log n \leq nR^2 < n^{2\alpha+1}$. Therefore

$$f(nR^2) \geq f(\log n) > 1.$$

Since $D \leq f(nR^2) + 1 \leq 2f(nR^2)$ and $D \geq f(nR^2)$, the right-hand term of (22) is thus bounded by

$$C\left\{ nR^2 (f(nR^2))^{-2\alpha} + 2f(nR^2)\left(1 + \log\left(\frac{n}{f(nR^2)}\right)\right) \right\}$$

and the result easily follows. □

**Proof of Theorem 5.** Let $s$ be an arbitrary element of $\mathcal{P}(R, \alpha, p)$ and recall that $|m| \leq c_0 2^J$ for every $m \in \mathcal{M}_J$, where $c_0 = 1 + \sum_{k \geq 1} k^{-3}$. If $c$ is large enough, then (2) holds with the $L_m$'s all equal to an absolute constant $L$, and $\Sigma \leq 1$ provided $L$ is large enough, see Section 5. It thus follows from Theorem 1, Lemma 5 and the definition of the penalty function pen that there exists $C > 0$ that only depends on $\alpha$, $p$ and $c$ such that for every $J = 0, \ldots, N - 1$,

$$\mathbb{E}_s \|s - \tilde{s}\|^2 \leq C\{nR^2 2^{-2\alpha J} + 2^J\}.$$

In order to optimize this inequality, we search for the $J$ that provides the smallest upper bound, that is, some $J$ that approximately equates the terms $nR^2 2^{-2\alpha J}$ and $2^J$: We consider

$$J = \sup\{j \in \mathbb{N} \text{ s.t. } 2^j \leq (nR^2)^{1/(2\alpha+1)}\}.$$

Then $J$ is well defined, since $nR^2 \geq 1$. Moreover, $J \leq N - 1$, since $n > (nR^2)^{1/(2\alpha+1)}$ (recall $n^\alpha > R$). The result easily follows, since

$$2^J \leq (nR^2)^{1/(2\alpha+1)} \quad \text{and} \quad 2^{-J-1} \leq (nR^2)^{-1/(2\alpha+1)}. \qquad \square$$

# Appendix

## A.1. Proof of Lemma 1

Let $S = S_m + S_{m'}$ for some $m' \in \mathcal{M}$ and $\varepsilon'$ be a $r \times n$ random matrix that has the same distribution as $\varepsilon$ and is independent of $\varepsilon$. Set

$$Z = \sup_{u \in \mathbb{R}^r \otimes S} \frac{\langle \varepsilon, u \rangle}{\|u\|} \tag{A.1}$$

and for every $i = 1, \ldots, n$ let us define $Z^{(i)}$ in the same way as $Z$, but with $\varepsilon_i$ replaced by $\varepsilon'_i$:

$$Z^{(i)} = \sup_{u \in \mathbb{R}^r \otimes S} \frac{1}{\|u\|} \left( \sum_{k \neq i} \langle \varepsilon_k, u_k \rangle_r + \langle \varepsilon'_i, u_i \rangle_r \right),$$

where $\langle \cdot, \cdot \rangle_r$ denotes the Euclidean inner product in $\mathbb{R}^r$. By the Cauchy–Schwarz inequality the supremum in (A.1) is achieved at an $\varepsilon$-measurable point $u^*$ (which is the orthogonal projection of $\varepsilon$ onto $\mathbb{R}^r \otimes S$). Hence

$$Z - Z^{(i)} \leq \frac{1}{\|u^*\|} \langle \varepsilon_i - \varepsilon'_i, u_i^* \rangle_r.$$

It is assumed that $\varepsilon'$ has the same distribution as $\varepsilon$ and is independent of $\varepsilon$. Since $\varepsilon'$ is centered and $u^*$ is $\varepsilon$-measurable, we get

$$\mathbb{E}_s[(\langle \varepsilon_i - \varepsilon'_i, u_i^* \rangle_r)^2 | \varepsilon] = \sum_{j=1}^r \sum_{j'=1}^r u_i^{*(j)} u_i^{*(j')} (\varepsilon_i^{(j)} \varepsilon_i^{(j')} + \mathbb{E}_s(\varepsilon_i^{(j)} \varepsilon_i^{(j')}))$$

$$= \sum_{j=1}^r (u_i^{*(j)})^2 (X_i^{(j)} - 2 X_i^{(j)} s_i^{(j)} + s_i^{(j)}) - 2 \sum_{\substack{j=1 \\ j \neq j(i)}}^r u_i^{*(j(i))} u_i^{*(j)} s_i^{(j)},$$

where $j(i)$ is the index such that $X_i^{(j(i))} = 1$ and $X_i^{(j)} = 0$ for all $j \neq j(i)$. On the one hand, distinguishing the cases $X_i^{(j)} = 0$ and $X_i^{(j)} = 1$, one has $X_i^{(j)} - 2X_i^{(j)} s_i^{(j)} + s_i^{(j)} \leq 1$. On the other hand, using the inequality $2\sqrt{ab} \leq a + b$ for all positive numbers $a$ and $b$ and Jensen's inequality, we get

$$\left| 2 \sum_{\substack{j=1 \\ j \neq j(i)}}^r u_i^{*(j(i))} u_i^{*(j)} s_i^{(j)} \right| \leq (u_i^{*(j(i))})^2 + \left( \sum_{\substack{j=1 \\ j \neq j(i)}}^r u_i^{*(j)} s_i^{*(j)} \right)^2 \leq \sum_{j=1}^r (u_i^{*(j)})^2.$$

Hence

$$\mathbb{E}_s \left[ \sum_{i=1}^n (Z - Z^{(i)})^2 \mathbb{1}_{Z > Z^{(i)}} | \varepsilon \right] \leq 2.$$

By Corollary 3 of [7] we thus have

$$\mathbb{P}_s\left(\sup_{u\in\mathbb{R}^r\otimes S}\frac{\langle\varepsilon,u\rangle}{\|u\|}>\mathbb{E}_s\left[\sup_{u\in\mathbb{R}^r\otimes S}\frac{\langle\varepsilon,u\rangle}{\|u\|}\right]+2\sqrt{2x}\right)\leq\exp(-x)$$

for all $x\geq 0$. Let $D$ be the dimension of $S$ and $\{\phi_1,\ldots,\phi_D\}$ be an orthonormal basis of $S$. Recall that the above supremum is the norm of the orthogonal projection of $\varepsilon$ onto $\mathbb{R}^r\otimes S$ so by Jensen's inequality,

$$\mathbb{E}_s\left[\sup_{u\in\mathbb{R}^r\otimes S}\frac{\langle\varepsilon,u\rangle}{\|u\|}\right]=\mathbb{E}_s\left[\left(\sum_{j=1}^r\sum_{k=1}^D(\langle\phi_k,\varepsilon^{(j)}\rangle_n)^2\right)^{1/2}\right]$$

$$\leq\left(\sum_{j=1}^r\sum_{k=1}^D\mathbb{E}_s(\langle\phi_k,\varepsilon^{(j)}\rangle_n)^2\right)^{1/2}.$$

But $D\leq D_m+D_{m'}$, $\mathrm{var}_s(\varepsilon_i^{(j)})\leq s_i^{(j)}$ for every $i$ and $j$ and $\|\phi_k\|_n^2=1$ for every $k$, so

$$\mathbb{E}_s\left[\sup_{u\in\mathbb{R}^r\otimes S}\frac{\langle\varepsilon,u\rangle}{\|u\|}\right]\leq\sqrt{D_m+D_{m'}}$$

and the result follows.

### A.2. Proof of Lemma 2

Let $\mu$ be the counting measure on $\mathbb{N}$. We consider the distance $\delta$ on the set of all subsets of $\{1,\ldots,n\}$ given by

$$\delta(m,m')=\tfrac{1}{2}\int|\mathbb{1}_m(x)-\mathbb{1}_{m'}(x)|\,\mathrm{d}\mu(x).$$

We begin by proving the following lemma.

**Lemma 6.** *Let $D$ and $D_0$ be positive integers with $D\leq D_0$. Assume that either $D_0\geq 18D$ or $D_0<18D$ and $D\geq 534$. Then there exists a set $M$ of subsets of $\{2,\ldots,D_0+1\}$ such that $|m|\leq D$ for all $m\in M$, $\delta(m,m')\geq D/7$ for all $m\neq m'\in M$, $|M|\geq 6$ and*

$$\log|M|\geq\frac{D}{150}\log\left(\frac{D_0}{D}\vee 18\right).$$

**Proof.** Assume first $D_0\geq 18D$. The result is trivial for $D=1$ (consider $M=\{\{j\},j=2,\ldots,D_0\}$) so we assume in the sequel $D\geq 2$. By Lemma 4 of [5], there exists a set $M$ of subsets of $\{2,\ldots,1+D_0\}$ such that $|m|=2\lfloor D/2\rfloor$ for all $m\in M$, $\delta(m,m')\geq D/2$ for all $m \neq m' \in M$ and

$$\log |M| \geq \lfloor D/2 \rfloor \left( \log \left( \frac{D_0}{\lfloor D/2 \rfloor} \right) - \log 16 + 1 \right).$$

Since $\lfloor D/2 \rfloor \geq 1$ and $D_0 \geq 18D$, we thus have $|M| \geq 6$. Moreover, $\log 8 - 1 \leq (\log(D_0/D))/2$ and $\lfloor D/2 \rfloor \geq D/3$, so

$$\log |M| \geq \frac{D}{6} \log \left( \frac{D_0}{D} \right).$$

This set $M$ thus satisfies all the required properties and the proof is complete in the case $D_0 \geq 18D$. Assume now $D_0 < 18D$ and $D \geq 534$. Let $k = \lfloor D/6 \rfloor$. We use Lemma 9 of [3] with $\Omega = \{2, \ldots, 6k+1\}$, $C = 3k$ and $q = k$. We know from this lemma that there exists a set $M$ of subsets of $\{2, \ldots, 6k+1\}$ such that $|m| = 3k$ for all $m \in M$, $\delta(m, m') \geq k$ for all $m \neq m' \in M$ and

$$|M| \geq \tfrac{3}{4} \binom{6k}{3k} \left[ \binom{3k}{2k} \binom{3k}{k} \right]^{-1}.$$

On the one hand,

$$(6k)! = \left[ \prod_{a=0}^{3k-1} (6k - 2a) \right] \times \left[ \prod_{a=0}^{3k-2} (6k - 2a - 1) \right] \geq [2^{3k}(3k)!] \times [2^{3k-1}(3k-1)!].$$

On the other hand, $(3k)! \leq (3^k k!)^3$ and $(2k)! \geq 2^{2k-1} k!(k-1)!$, so

$$|M| \geq \frac{1}{32} \times \frac{2^{10k}}{k^3 3^{6k}}.$$

Also, $k^3 \leq 10^4 (2^8/3^5)^k$ and $32 \times 10^4 \leq (4/3)^{k/2}$ since $k \geq 89$ (recall $D \geq 534$), so $|M| \geq (4/3)^{k/2}$. It is then easy to check that $M$ satisfies all the required properties. $\square$

We turn now to the proof of Lemma 2, distinguishing between two cases.

Case 1. Either $D_0 \geq 18D$ or $D_0 < 18D$ and $D \geq 534$. Let $A > 0$ and $D_0, D$ integers with $1 \leq D \leq D_0 < n$. Assume that either $D_0 \geq 18D$ or $D_0 < 18D$ and $D \geq 534$. For every $m \subset \{2, \ldots, D_0 + 1\}$ let $s_m$ be the $r \times n$ matrix defined by

$$\begin{cases} s_m^{(1)} = 1/2 + a \sum_{\lambda \in m} \phi_\lambda, \\ s_m^{(2)} = 1/2 - a \sum_{\lambda \in m} \phi_\lambda, \\ s_m^{(j)} = 0, & \text{for } j > 2 \text{ (in the case } r > 2\text{),} \end{cases}$$

where

$$a = \frac{A}{\sqrt{2}} \wedge \left( \frac{1}{60} \sqrt{\log \left( \frac{D_0}{D} \vee 18 \right)} \right) \wedge \left( \frac{1}{20} \sqrt{\frac{n}{D_0}} \right).$$

Let $\|\cdot\|_\infty$ denote the supremum norm in $\mathbb{R}^n$. By definition, $\|\phi_\lambda\|_\infty \leq 2^{j/2}/\sqrt{n}$ for every $\lambda \in \Lambda(j)$ and $j = 0, \ldots, N-1$. Moreover, the functions $\phi_\lambda$, $\lambda \in \Lambda(j)$ have disjoint supports so

$$\left\|\sum_{\lambda \in \Lambda(j)} \phi_\lambda\right\|_\infty \leq \frac{2^{j/2}}{\sqrt{n}}$$

for every $j$. Let $J_0$ be that integer such that $2^{J_0} < D_0 + 1 \leq 2^{J_0+1}$. Then,

$$\{2, \ldots, D_0 + 1\} \subset \bigcup_{j=0}^{J_0} \Lambda(j)$$

(recall $|\Lambda(j)| = 2^j$ for $j \neq -1$) so for every $m \subset \{2, \ldots, 1 + D_0\}$,

$$\left\|\sum_{\lambda \in m} \phi_\lambda\right\|_\infty \leq \sum_{j=0}^{J_0} \frac{2^{j/2}}{\sqrt{n}} \leq 5\sqrt{\frac{D_0}{n}}.$$

Therefore, $a\|\sum_{\lambda \in m} \phi_\lambda\|_\infty \leq 1/4$ for all $m \subset \{2, \ldots, D_0 + 1\}$. It follows that $s_m \in \tilde{\mathcal{P}}(A, D_0, D)$ for every subset $m$ of $\{2, \ldots, D_0 + 1\}$ with cardinality $|m| \leq D$. We now use Proposition 9 of [5], page 263 with $\mathcal{C} = \{s_m, m \in M\}$ and $M$ is as in Lemma 6. By Lemma 4,

$$K(P_{s_m}, P_{s_{m'}}) \leq 4\|s_m - s_{m'}\|^2$$

for every $m \neq m' \in M$, and since the vectors $\phi_\lambda$ are orthonormal, we have

$$\|s_m - s_{m'}\|^2 = \sum_{j=1}^{2} \|s_m^{(j)} - s_{m'}^{(j)}\|^2 = 4a^2 \delta(m, m').$$

Therefore, $K(P_{s_m}, P_{s_{m'}}) \leq 16a^2 D$ and $\|s_m - s_{m'}\|^2 \geq 4a^2 D/7$ for every $m \neq m' \in M$. We thus set $\eta^2 = a^2 D/2$ and $H = 16a^2 D$. By definition of $a$,

$$\log|M| \geq \frac{D}{150} \log\left(\frac{D_0}{D} \vee 18\right) \geq \frac{D}{150} 3600a^2.$$

Therefore, $H/\log|M| \leq 2/3$ and Proposition 9 of [5] yields

$$\inf_{\hat{s}} \sup_{s \in \mathcal{C}} \mathbb{E}_s \|s - \hat{s}\|^2 \geq \frac{a^2 D}{24},$$

where the infimum is taken over all estimators $\hat{s}$ such that each column $\hat{s}^{(j)}$ is a probability on $\{1, \ldots, r\}$. Since $\mathcal{C} \subset \tilde{\mathcal{P}}(A, D_0, D)$, we then get (19) by extending the infimum over all possible estimators.

Case 2. $D_0 < 18D$ and $D \leq 533$. Let $A > 0$ and $D_0, D$ integers with $1 \leq D \leq D_0 < n$, $D_0 < 18D$ and $D \leq 533$. Let $v = s_m$ where $s_m$ is defined as in Case 1 with $m = \{2\}$ and

$$a = \frac{A}{\sqrt{2}} \wedge \frac{1}{4} \wedge \frac{\sqrt{n}}{4}.$$

We have $v \in \tilde{\mathcal{P}}(A, D_0, D)$ and, similar to Case 1, $\|u - v\|^2 = 2a^2$ and $K(P_u, P_v) \leq 8a^2$. By Lemma 3,

$$\inf_{\hat{s}} \sup_{s \in \{u,v\}} \mathbb{E}_s \|s - \hat{s}\|^2 \geq \frac{a^2}{2}(1 - \sqrt{8a^2}) \geq \frac{a^2}{4}$$

and the result follows from the assumptions $D_0 < 18D$ and $D \leq 533$.